\documentclass[12pt,leqno]{article}
\usepackage{amsfonts}
\usepackage{amsmath, amsthm, amsfonts, amssymb, color}
\usepackage{mathrsfs}
\usepackage{color}
\usepackage[pagewise]{lineno}\linenumbers
\theoremstyle{definition}
\newtheorem{defn}{Definition}[section]
\newcommand{\scr}[1]{\mathscr #1}
\definecolor{wco}{rgb}{0.5,0.2,0.3}

\numberwithin{equation}{section} \theoremstyle{remark}

\newcommand{\ua}{\uparrow}

\title{{\bf     Derivative Formula   for  Singular McKean-Vlasov  SDEs }\footnote{Supported in
 part by  NNSFC (11771326, 11831014, 11921001).} }
\author{{\bf  Feng-Yu Wang   }\\
\footnotesize{ Center for Applied Mathematics, Tianjin University, Tianjin 300072, China}\\
\footnotesize{    wangfy@tju.edu.cn}}
\begin{document}
\allowdisplaybreaks
\def\R{\mathbb R}  \def\ff{\frac} \def\ss{\sqrt} \def\B{\mathbf
B}
\def\N{\mathbb N} \def\kk{\kappa} \def\m{{\bf m}}
\def\ee{\varepsilon}\def\ddd{D^*}
\def\dd{\delta} \def\DD{\Delta} \def\vv{\varepsilon} \def\rr{\rho}
\def\<{\langle} \def\>{\rangle}
  \def\nn{\nabla} \def\pp{\partial} \def\E{\mathbb E}
\def\d{\text{\rm{d}}} \def\bb{\beta} \def\aa{\alpha} \def\D{\scr D}
  \def\si{\sigma} \def\ess{\text{\rm{ess}}}\def\s{{\bf s}}
\def\beg{\begin} \def\beq{\begin{equation}}  \def\F{\scr F}
\def\Ric{\mathcal Ric} \def\Hess{\text{\rm{Hess}}}
\def\e{\text{\rm{e}}} \def\ua{\underline a} \def\OO{\Omega}  \def\oo{\omega}
 \def\tt{\tilde}\def\[{\lfloor} \def\]{\rfloor}
\def\cut{\text{\rm{cut}}} \def\P{\mathbb P} \def\ifn{I_n(f^{\bigotimes n})}
\def\C{\scr C}      \def\aaa{\mathbf{r}}     \def\r{r}
\def\gap{\text{\rm{gap}}} \def\prr{\pi_{{\bf m},\varrho}}  \def\r{\mathbf r}
\def\Z{\mathbb Z} \def\vrr{\varrho} \def\ll{\lambda}
\def\L{\scr L}\def\Tt{\tt} \def\TT{\tt}\def\II{\mathbb I}
\def\i{{\rm in}}\def\Sect{{\rm Sect}}  \def\H{\mathbb H}
\def\M{\mathbb M}\def\Q{\mathbb Q} \def\texto{\text{o}} \def\LL{\Lambda}
\def\Rank{{\rm Rank}} \def\B{\scr B} \def\i{{\rm i}} \def\HR{\hat{\R}^d}
\def\to{\rightarrow}\def\l{\ell}\def\iint{\int}\def\gg{\gamma}
\def\EE{\scr E} \def\W{\mathbb W}
\def\A{\scr A} \def\Lip{{\rm Lip}}\def\S{\mathbb S}
\def\BB{\scr B}\def\Ent{{\rm Ent}} \def\i{{\rm i}}\def\itparallel{{\it\parallel}}
\def\g{{\mathbf g}}\def\Sect{{\mathcal Sec}}\def\T{\mathcal T}\def\BB{{\bf B}}
\def\f{\mathbf f} \def\g{\mathbf g}\def\BL{{\bf L}}  \def\BG{{\mathbb G}}
\def\Bd{{D^E}} \def\BdP{D^E_\phi} \def\Bdd{{\bf \dd}} \def\Bs{{\bf s}} \def\GA{\scr A}
\def\Bg{{\bf g}}  \def\Bdd{\psi_B} \def\supp{{\rm supp}}\def\div{{\rm div}}
\def\ddiv{{\rm div}}\def\osc{{\bf osc}}\def\1{{\bf 1}}\def\BD{\mathbb D}\def\GG{\Gamma}
\def\H{{\bf H}}
\maketitle

\begin{abstract}   The Bismut formula is established   for the intrinsic derivative of singular  McKean-Vlasov SDEs, where  the noise coefficient belongs to a local Sobolev space,
 and the  drift contains a  locally integrable time-space term as well as  a  time-space-distribution term Lipschitz continuous in the space and distribution variables.  The results are   new  also for   classical SDEs.
  \end{abstract} \noindent

 \vskip 2cm

 \section{Introduction and main results}

 Since 1984 when Bismut \cite{Bismut} presented  his derivative formula for diffusion semigroups on Riemannian manifolds, this type formula has been widely developed and applied. Recently, Bismut formula was established in \cite{XXZZ}
 for singular SDEs with a locally integrable drift.

 On the other hand, as crucial probability models characterizing nonlinear
 Fokker-Planck equations and mean field games, distribution dependent (also called McKean-Vlasov or mean-field) SDEs have been intensively investigated, see for instance the monographs \cite{SN, CD} and the survey \cite{20HRW}.
  In particular, Bismut type formulas
 have been established  in \cite{RW19, BRW20, HSW21} for
regular  McKean-Vlasov SDEs where the drift is at least Dini continuous in the space variable,  see also \cite{DM} for Bismut formula of the decoupled SDE where the distribution parameter is fixed.

In this paper, we aim to establish  Bismut formula of the intrinsic  derivative in the $L^k$-Wasserstein space ($k\ge 1)$  for singular McKean-Vlasov SDEs. In the following we first introduce
the model considered in the paper, then  recall the intrinsic   derivative, and finally state the main results of the paper.

\subsection{McKean-Vlasov SDE}

Let $\scr P$ be the space of probability measures on $\R^d$.  We will use $|\cdot|$ to denote the absolute value or the norm in the Euclidean space, and  $\|\cdot\|$  the operator norm for  linear operators or matrices.
 For any $k\in [1,\infty)$, the $L^k$-Wasserstein space
$$\scr P_k:=\{\mu\in \scr P:   \mu(|\cdot|^k)<\infty\}$$ is a Polish space
under  the  $L^k$-Wasserstein distance
$$\W_k(\mu,\nu):=   \inf_{\pi\in \C(\mu,\nu)} \bigg(\int_{\R^d\times\R^d} |x-y|^k \pi(\d x,\d y) \bigg)^{\ff 1{k}},\ \ \mu,\nu\in \scr P_ k,$$
where $\C(\mu,\nu)$ is the set of all couplings for $\mu$ and $\nu$.

Throughout the paper, we fix   $T\in (0,\infty)$ and  consider the following McKean-Vlasov SDE on $\R^d$:
 \beq\label{E10} \d X_t=  b_t(X_t,\L_{X_t}) \d t+\si_t(X_t)\d W_t,\ \ t\in [0,T],\end{equation}
where $W_t$ is an $m$-dimensional Brownian motion on a complete filtered  probability space $(\OO, \{\F_t\}_{t\ge 0},\P)$,   $\L_\xi$ is the distribution (i.e. the law) of a random variable $\xi$,
and \beg{align*}b: [0,\infty)\times\R^d\times \scr  P_k\to \R^d,
 \ \ \si: [0,\infty)\times\R^d \to \R^d\otimes\R^m\end{align*}  are measurable. When different probability measures are considered, we denote $\L_{\xi}=\L_{\xi|\P}$ to emphasize the distribution of $\xi$ under $\P$.

\beg{defn}  (1) A continuous adapted process  $(X_{t})_{t\in [0,T]}$ is called a  solution of \eqref{E10}, if  with $\L_{X_\cdot}\in C([0,T];\scr P_k)$,
$$\int_0^T\E\big[|b_r(X_{r}, \L_{X_{r}})|+  \| \si_r(X_{r})\|^2\big]\d r<\infty $$ and   $\P$-a.s.
  $$X_{t} = X_{0} +\int_0^t b_r(X_{r}, \L_{X_{r}})\d r + \int_0^t \si_r(X_{r})\d W_r,\ \ t\in [0,T].$$
 We call $\eqref{E10}$ strongly well-posed for distributions in $\scr P_k$, if it has a unique solution for any initial value $X_0\in L^k(\OO\to\R^d,\F_0,\P)$.

 (2) A couple $(\tt X_{t},\tt W_t)_{t\in [0,T]}$   is called a  weak solution of \eqref{E10}, if   $\tt W_t$  is an   $m$-dimensional Brownian motion on a complete filtered probability space $(\tt\OO,\{\tt\F_t\}_{t\in [0,T]}, \tt\P)$  such that   $(\tt X_{t})_{t\in [0,T]}$  is a solution of \eqref{E10} for   $(\tt W_t, \tt\P)$ replacing   $(W_t,\P)$.    \eqref{E10} is called weakly well-posed for distributions in $\scr P_k$, if  any   $\nu\in\scr  P_k$
 it has a weak solution with initial distribution $\nu$, and for any two weak solutions $(X_t^i, W_t^i)$ under $\P^i$, $ i=1,2$, $\L_{X_0^1|\P^1}=\L_{X_0^2|\P^2}$ implies $\L_{(X^1_t)_{t\in [0,T]}|\P^1}= \L_{(X^2_t)_{t\in [0,T]}|\P^2}$.

(3) We call  \eqref{E10} well-posed for distributions in $\scr P_k$, if it is both strongly and weakly well-posed for distributions in $\scr P_k$.
 \end{defn}

When the SDE \eqref{E10} is well-posed   for distributions in $\scr P_k$,    we denote  $P_t^*\mu=\L_{X_t}$ for the solution with initial distribution $\mu\in \scr P_k$, and define the family of linear operators $\{P_t\}_{t\in [0,T]}$ from $\B_b(\R^d)$ to $\B_b(\scr P_k)$, where $\B_b(\cdot)$ stands for the set of all bounded measurable functions on a measurable space:
$$P_t f(\mu):= (P_t^*\mu)(f)= \int_{\R^d} f\d(P_t^*\mu),\ \ t\in [0,T], f\in \B_b(\R^d), \mu\in \scr P_k.$$

To characterize  the singularity of coefficients $b$ and $\si$  in time-space variables, we recall some functional spaces introduced in \cite{XXZZ}.  For any $p\ge 1$, $L^p(\R^d)$ is the class of   measurable  functions $f$ on $\R^d$ such that \index{$L^p(\R^d)$: $L^p$-space for Lebesgue measure on $\R^d$}
 $$\|f\|_{L^p(\R^d)}:=\bigg(\int_{\R^d}|f(x)|^p\d x\bigg)^{\ff 1 p}<\infty.$$
For any $p,q>1$, let $\tt L_q^p$ denote the class of measurable functions $f$ on $[0,T]\times\R^d$ such that
$$\|f\|_{\tt L_q^p}:= \sup_{z\in \R^d}\bigg( \int_{0}^{T} \|1_{B(z,1)}f_t\|_{L^p(\R^d)}^q\d t\bigg)^{\ff 1 q}<\infty,$$
where $B(z,1):= \{x\in\R^d: |x-z|\le 1\}$. We denote $f\in \tt H_{q}^{2,p}$ if $|f|+|\nn f|+|\nn^2f|\in \tt L_{q}^p$.
We will  take $(p,q)$ from the class
$$\scr K:=\Big\{(p,q): p,q\in  (2,\infty),\  \ff d p+\ff 2 q<1\Big\}.$$ \index{$\scr K:=\{(p,q): p,q\in (1,\infty), \ff d p+\ff 2 q<1\}$}

\subsection{Intrinsic derivative in $\scr P_k$}

The intrinsic derivative for measures was introduced in \cite{AKR} to  construct diffusion processes on configuration spaces over   Riemannian manifolds,    and   used in \cite{OTTO} to study the geometry of dissipative evolution equations, see \cite{AGS} for analysis and geometry on the Wasserstein space over a metric measure space. This derivative corresponds to the motion of particles, comparing to the extrinsic derivative induced by the birth and death of particles.
  The following   notion   was introduced in \cite{BRW20} for functions  on the $L^k$-Wasserstein space over a Banach space.

   For any $\mu\in \scr P_k$, the tangent space at $\mu$ is
 $$T_{\mu,k} := L^k(\R^d\to\R^d;\mu).$$
 Then $\mu\circ(id+\phi)^{-1} \in \scr P_k$ for $\phi\in T_{\mu,k}$, where $id$ is the identity map on $\R^d$.

\begin{defn} Let  $f$  be a continuous function  on $\scr P_k$.
It  is called intrinsically differentiable at a point $\mu\in\scr P_k$, if
 $$T_{\mu,k}\ni\phi\mapsto D_\phi^If(\mu):= \lim_{\vv\downarrow 0} \ff{f(\mu\circ(id+\vv\phi)^{-1})-f(\mu)}{\vv}\in\R
$$ is a well  defined bounded linear functional. In this case,  the intrinsic derivative is the unique element
$$D^If(\mu)\in T_{\mu,k}^*:=L^{k^*}(\R^d\to\R^d;\mu),\ \ k^*:=\ff k{k-1}\big(=\infty\ \text{if}\ k=1\big)$$ such that \index{$(D^If(\mu)$: intrinsic derivative of $f$ at $\mu$}
$$\int_{\R^d}\<D^If(\mu)(x), \phi(x)\>  \mu(\d x) = D_\phi^I f(\mu),\ \ \phi\in T_{\mu,k}.$$
If moreover
 \begin{equation*}
\lim_{\|\phi\|_{T_{\mu,k}}\downarrow0}\ff{|f(\mu\circ(id+\phi)^{-1})-f(\mu)-D_\phi^If(\mu)|}{\|\phi\|_{T_{\mu,k}}}=0,
\end{equation*}
then $f$ is called    $L$-differentiable at $\mu$.
The function  $f$ is called intrinsically (or $L$-) differentiable on $\scr P_k$, if it is intrinsically (or $L$-) differentiable at any $\mu\in \scr P_k$.
 \end{defn}

 Note that when $k=2$ the Lions derivative $D^Lf(\mu)$ is defined as the unique element in $T_{\mu,2}$ such that for any atomless probability space $(\OO,\F,\P)$ and any random variables $X,Y$ with $\L_{X}=\mu$,
 $$\lim_{\|Y-X\|_{L^2(\P)}\downarrow 0}\ff{|f(\L_Y)-f(\L_X)-\E[\<D^Lf(\mu)(X),Y-X\>]|}{\|Y-X\|_{L^2(\P)}}=0.$$
 Since $D^Lf(\mu)$ does not depend on the choice of probability space, when $\mu$ is atomless we may choose
 $(\OO,\F,\P)=(\R^d,\B^d,\mu)$ such that $D^Lf(\mu)=D^If(\mu)$, see for instance \cite[Chapter 5]{CD}.  So, in the following, when $f$ is $L$-differentiable we denote $D^If=D^Lf$.

To measure the singularity of the drift,  we write  $b_t(x,\mu)= b_t^{(0)}(x)+ b_t^{(1)}(x,\mu)$, where $|b^{(0)}|\in \tt L_{q_0}^{p_0}$ for some $(p_0,q_0)\in \scr K$ and $b_t^{(1)}$ is in the class $\D_k$ defined
as follows.

\beg{defn}   $\D_k$ is the class  of continuous functions $g$ on $\R^d\times\scr P_k$ such that $g(x,\mu)$ is differentiable in $x$, $L$-differentiable in $\mu$, and
$D^Lg(x,\mu)(y)$ has a version jointly continuous in $(x,y,\mu)\in \R^d\times \R^d\times\scr P_k$ such that
$$|D^Lg(x,\mu)(y)|\le c(x,\mu) (1+|y|^{k-1}),\ \ x,y\in\R^d,\mu\in \scr P_k$$ holds for some positive function $c$ on $\R^d\times\scr P_k.$  \end{defn}

Typical examples of functions in $\D_k$ are cylindrical functions of type
$$g(x,\mu):=F(x,\mu(h_1),\cdots, \mu(h_n))$$
for some $F\in C^1(\R^d\times\R^n)$ and $\{h_i\}_{1\le i\le n}\subset C^1(\R^d)$ such that
$$\sup_{1\le i\le n} |\nn h_i(y)|\le c (1+|y|^{k-1}),\ \ y\in \R^d$$ holds for some constant $c>0$. In this case
$$ D^L g(x,\mu)(y) =\sum_{i=1}^n\{\pp_i F(x,\cdot)\}(\mu(h_1),\cdots, \mu(h_n)) \nn h_i(y),\ \ x,y\in\R^d, \mu\in \scr P_k.$$

 \subsection{Main results}

  To establish Bismut  formula of $P_tf$  on $\scr P_k$, we make the following assumption.

 \beg{enumerate}\item[$(H)$]      $b_t(x,\mu)=b^{(0)}_t(x)+ b_t^{(1)}(x,\mu) $ such that the following conditions hold.
  \item[$(1)$]  $a:= \si\si^*$ is invertible with $\|a\|_\infty+\|a^{-1}\|_\infty<\infty$, where $\si^*$ is the transposition of $\si$,   and
$$ \lim_{\vv\to 0} \sup_{|x-y|\le \vv, t\in [0,T]} \|a_t(x)-a_t(y)\|=0.$$
 \item[$(2)$]   $|b^{(0)}|\in \tt L_{q_0}^{p_0}$ for some  $(p_0,q_0)\in \scr K$.  Moreover, $\si_t$ is a.e. differentiable such that
\beq\label{ASI} \|\nn \si\|\le \sum_{i=1}^l f_i \end{equation}  holds for some $l\in \mathbb N$ and $0\le f_i \in \tt L_{q_i}^{p_i}$ with $  (p_i,q_i) \in  \scr K, 1\le i\le l.$
 \item[$(3)$]          $b_t^{(1)} \in \D_k$ such that
 \beq\label{LPS} \sup_{(t,x,\mu)\in [0,T]\times\R^d\times \scr P_k} \big\{|b_t^{(1)}(0,\dd_0)|+ \|\nn  b_t^{(1)}(x,\mu)\| + \|D^L b_t^{(1)}(x,\mu) \|_{L^{k^*}(\mu)} \big\}<\infty, \end{equation}
 where   $\dd_0$ is the Dirac measure at $0\in\R^d$, $\nn$ is the gradient in the space variable $x\in\R^d$, and
  $D^L$ is the $L$-derivative in the distribution variable $\mu\in \scr P_k$.
\end{enumerate}

We will show that $(H)$   implies the well-posedness of \eqref{E10} for distributions in $\scr P_k$.
To calculate the intrinsic derivative $D^I P_tf(\mu)$, for any $\vv\in [0,1]$ and $\phi\in T_{\mu,k}$, we consider the following  SDE:
 \beq\label{E101}   \beg{split} \d X_t^{\mu,\vv\phi}= &\  b_t(X_t^{\mu,\vv\phi},\L_{X_t^{\mu,\vv\phi}}) \d t+\si_t(X_t^{\mu, \vv\phi})\d W_t,\\
 & \  t\in [0,T],\ \ \ X_0^{\mu,\vv\phi}= X_0^\mu+ \vv \phi(X_0^\mu).\end{split} \end{equation}
 Note that $X_t^\mu:=X_t^{\mu,0}$ solves \eqref{E10}.
 We will prove that the derivative process
 \beq\label{DPHI} \nn_\phi X_t^\mu:= \lim_{\vv\downarrow 0} \ff{X_t^{\mu,\vv\phi} -X_t^\mu}\vv,\ \ t\in [0,T]\end{equation}
 exists in $L^k(\OO\to C([0,T];\R^d),\P)$. We also need the derivative of the decoupled SDE
  \beq\label{E102}  \beg{split} & \d   X_t^{\mu,x}=  b_t(X_t^{\mu,x},P_t^*\mu) \d t+\si_t(X_t^{\mu,x})\d W_t,\\
  & \   t\in [0,T],  X_0^{\mu,x}= x, x\in\R^d, \mu\in \scr P_k.\end{split}\end{equation}
 By Theorem \ref{T03},  $(H)$  implies the well-posedness of \eqref{E102} and that
 for any $v\in \R^d$,
 \beq\label{DPHI'} \nn_v X_t^{\mu,x}:= \lim_{\vv\downarrow 0} \ff{X_t^{\mu,x+\vv v(x)} -X_t^{\mu,x}}\vv,\ \ t\in [0,T]\end{equation}
 exists in $L^k(\OO\to C([0,T];\R^d),\P)$.

   Our first result present some continuity estimates in terms of the initial data.

\beg{thm}\label{TA2} Assume    $(H)$.  Then the following assertions hold.
 \beg{enumerate}  \item[$(1)$]
$\eqref{E10}$ is well-posed for distributions in $\scr P_k$, and for any $j\ge 1$ there exists a constant $c>0$ such that any solution $X_t$ satisfies
 \beq\label{EST0} \E \bigg[\sup_{t\in [0,T]} |X_t|^j\bigg|\F_0\bigg]\le c\Big\{1+ (\E[|X_0|^k])^{\ff j k}+   |X_0|^j\Big\}.\end{equation}  In particular, there exists a constant $c>0$ such that
 \beq\label{*BN} \E \bigg[\sup_{t\in [0,T]} |X_t|^k \bigg]\le c \big(1+ \E[|X_0|^k]\big).\end{equation}
\item[$(2)$]  For any $j\ge 1$ there exists a constant $c>0$ such that for any two solutions $X_t^1, X_t^2$ of $\eqref{E10}$
with  initial distributions in $\scr P_k$,
\beq\label{EE0} \E \bigg[\sup_{t\in [0,T]} |X_t^1-X_t^2|^j\bigg|\F_0\bigg] \le c\Big\{\big(\E[|X_0^1-X_0^2|^k]\big)^{\ff j k}+|X_0^1-X_0^2|^j\Big\}.\end{equation}
In particular, there exists a constant $c>0$ such that
 \beq\label{EE0*}\E \bigg[\sup_{t\in [0,T]} |X_t^1-X_t^2|^k\bigg] \le c   \E[|X_0^1-X_0^2|^k].\end{equation}
 \item[$(3)$] There exists a constant $c>0$ such that
 \beq\label{GRD00} \beg{split} &\|P_t^*\mu- P_t^*\nu\|_{var} := \sup_{\|f\|_\infty\le 1} |P_tf(\mu)-P_tf(\nu)|\\
 &\le \ff c {\ss t} \W_k(\mu,\nu),\ \ t\in (0,T], \mu,\nu\in \scr P_k.\end{split}\end{equation}
 \end{enumerate}
\end{thm}

By taking $X_0^1$ and $X_0^2$ such that
$$\L_{X_0^1}=\mu, \ \ \L_{X_0^2}=\nu,\ \ \E[|X_0^1-X_0^2|^k]=\W_k(\mu,\nu)^k,$$ we deduce from \eqref{EE0*} that
\beq\label{LN2}   \W_k(P_t^* \mu, P_t^*\nu)\le c \W_k(\mu,\nu),\ \ t\in [0,T], \mu,\nu\in \scr P_k\end{equation}  holds for some constant $c>0$.
Our next result provides derivative  estimates and Bismut formula of   $D^I P_t f$ for $t\in (0,T]$ and $f\in \B_b(\R^d)$.

\beg{thm}\label{TA2'} Assume    $(H)$.  Then the following assertions hold.
 \beg{enumerate}  \item[$(1)$]    For any $\mu \in \scr P_k$,  $\phi\in T_{\mu,k}$ and $v,x\in \R^d$,  $\nn_\phi X_t^\mu$ and $\nn_v X_t^{\mu,x}$ exist in $L^k(\OO\to C([0,T];\R^d),\P).$
Moreover,     for any $j\ge 1$ there exists a   constant $c>0$ such that
\beq\label{DFF}  \E\Big[\sup_{t\in [0,T]}|\nn_\phi X_t^\mu|^j\Big|\F_0\Big]\le c \big\{\|\phi\|_{L^k(\mu)}^j+ |\phi(X_0^\mu)|^j\big\},\ \ \mu\in \scr P_k, \phi\in T_{\mu,k},\end{equation}
\beq\label{DFF'}  \E\Big[\sup_{t\in [0,T]}|\nn_v X_t^{\mu,x}|^j\Big]\le c |v|^j,\ \ \mu\in \scr P_k,  x,v\in\R^d.\end{equation}
\item[$(2)$]  Denote $\zeta=\si(\si\si^*)^{-1}$. For any $t\in (0,T]$ and $f\in \B_b(\R^d)$, $P_t f$ is intrinsically differentiable on $\scr P_k$. Moreover, for any $\phi\in T_{\mu,k}$ and  $\bb \in C^1([0,t])$ with $\bb_0=0$ and $\bb_t=1$,
\beq\label{BSMI} \beg{split} & D_\phi^I  P_tf(\mu)= \int_{\R^d} \E\bigg[f(X_t^{\mu,x}) \int_0^t \bb_s'\big\<\zeta_s(X_s^{\mu,x} ) \nn_{\phi (x)}X_s^{\mu,x},  \d W_s\big\> \bigg]\mu(\d x)\\
&+  \E\bigg[f(X_t^{\mu}) \int_0^t  \Big\<\zeta_s(X_s^\mu)  \E\big[\<D^Lb^{(1)}_s(z, P_s^*\mu)(X_s^\mu), \nn_\phi X_s^\mu\>\big]|_{z=X_s^\mu},  \d W_s\Big\> \bigg].
\end{split} \end{equation}
Consequently, for any $p>1$   there exists a constant $c >0$ such that
\beq\label{GRDI}\beg{split} & \|D^I P_tf(\mu)\|_{L^{k^*}(\mu)} \le \ff {c}{\ss t} \Big\|\big(\E[|f|^p(X_t^\mu)|\F_0]\big)^{\ff 1 p} \Big\|_{L^{k^*}(\P)},\\
&  t\in (0,T], f\in \B_b(\R^d),\mu\in \scr P_k.\end{split} \end{equation}In particular, there exists a constant $c>0$ such that
$$\|D^I P_tf(\mu)\|_{L^{k^*}(\mu)} \le \ff {c}{\ss t} \big\|f(X_t^\mu)\big\|_{L^{k^*}(\P)},  \
  \ t\in (0,T], f\in \B_b(\R^d),\mu\in \scr P_k.$$
  \end{enumerate} \end{thm}

Finally, to prove the $L$-differentiability of $P_tf$, we need the uniform continuity of $\si_t(x)$, $\nn b_t^{(1)}(x,\mu)$ and $D^L b_t(x,\mu)(y)$ in $(x,y,\mu)$:
\beq\label{00*} \beg{split} &\lim_{\vv\downarrow 0}\sup_{t\in [0,T]} \sup_{|x-x'|\lor \W_k(\mu,\nu)\le \vv}\big\{\|\si_t(x)-\si_t(x')\|+ \big\|\nn b_t^{(1)}(x,\mu)- \nn b_t^{(1)}(x',\nu)\big\|\big\}=0,\\
&\lim_{\vv\downarrow 0}\sup_{t\in [0,T]} \sup_{|x-x'|\lor \W_k(\mu,\nu)\lor |y-y'|\le \vv} \big\|D^Lb_t^{(1)}(x,\mu)(y)- D^L b_t^{(1)}(x',\nu)(y')\big\|=0.\end{split} \end{equation}
Under this condition and $(H)$, the following result ensures the $L$-differentiability of $P_tf$ in $\scr P_k$ for $k>1$,  which improves the corresponding result in \cite{HSW21} where $k=2,$ $\si_t$ is Lipschitz continuous and $b^{(0)}_t$ is Dini continuous.

\beg{thm}\label{TA2''} Assume $(H)$ and $\eqref{00*}$ for $k\in (1,\infty)$. Then for any $t\in (0,T]$ and $f\in \B_b(\R^d)$, $P_tf$ is $L$-differentiable on $\scr P_k.$
\end{thm}

    In Section 2, we establish the Bismut  formula for singular SDEs where  $b_t(x,\mu)=b_t(x)$ is independent of $\mu$. Then we prove
   the above theorems  in Sections 3-5 respectively.

\section{Bismut formula for singular SDEs}

Let $b_t(x,\mu)=b_t(x)$ do not depend on $\mu$, so that \eqref{E10} becomes
 \beq\label{E10'} \d X_t= b_t(X_t)\d t+\si_t(X_t)\d W_t,\ \ t\in [0,T].\end{equation}
Let $X_t^x$ solves \eqref{E10'} for $X_0^x=x$, and consider
$$P_t f(x):= \E [f(X_t^x)],\ \ t\ge 0, f\in \B_b(\R^d), x\in \R^d.$$

Under $(H)$ with   $l=1$ and $b^{(1)}=0$,    the following Bismut formula \eqref{BSM}    is included in   Theorem 1.1(iii) of \cite{XXZZ} for   $f\in C_b^1(\R^d)$.
It is reasonable but nontrivial to extend   the formula from $f\in C_b^1(\R^d)$ to    $f\in \B_b(\R^d)$. The technique we used  in  step (d) in the proof of Theorem \ref{T03}(2) is due to \cite{RW19},
which will also be used in the proofs of Theorems \ref{TA2'} and \ref{TA2''}.

\beg{thm}\label{T03} Let  $(H)$ hold for $b_t^{(1)}(x,\mu)=b^{(1)}_t(x)$. Then $\eqref{E10'}$ is well-posed and the following assertions hold.
\beg{enumerate} \item[$(1)$] For any $j\ge 1$ and $x,v \in \R^d$,
$$\nn_vX_t^x:=\lim_{\vv\downarrow  0} \ff{X_t^{x+\vv v}-X_t^x}\vv,\ \ t\in [0,T]$$
exists in $L^j(\OO\to C([0,T];\R^d),\P)$, and  there exists a constant $c(j)>0$    such that
\beq\label{DFF0} \sup_{x\in \R^d} \E\Big[\sup_{t\in [0,T]}|\nn_v X_t^x|^j\Big]\le c(j) |v|^j,\ \  v\in\R^d.\end{equation}
\item[$(2)$]  For any $t\in (0,T]$ and $\bb\in C^1([0,t])$ with $\bb_0=0$ and $\bb_t=1$,
\beq\label{BSM} \nn_v P_tf(x)= \E\bigg[f(X_t^x) \int_0^t \bb_s'\Big\<\big\{\si_s^*(\si_s \si_s^*)^{-1} \big\}(X_s^{x}) \nn_vX_s^x,  \d W_s\Big\> \bigg]\end{equation}
holds for any $x\in \R^d$ and $f\in \B_b(\R^d).$
Consequently, for any $p>1$ there exists a constant $c(p)>0$ such that
\beq\label{GRD} |\nn P_tf|\le \ff {c(p)}{\ss t} \big(P_t |f|^p(x)\big)^{\ff 1 p},\ \ t\in (0,T], f\in \B_b(\R^d).\end{equation}\end{enumerate}
\end{thm}

\beg{proof}    The well-posedness follows from Lemma 3.1 in \cite{W21c}. Below we prove assertions (1) and (2) by using Zvonkin's transform as in \cite{XXZZ}.

(a) Zvonkin's transform.   Let  $\nn_v$ be the directional derivative along $v$. Consider
$$L_t=\ff 1 2{\rm tr}\big\{\si_t\si_t^*\nn^2\big\}+\nn_{b_t},\ \ t\in [0,T].$$
By \cite[Theorem 2.1]{YZ0}, there exists $\ll_0>0$, such that when $\ll\ge \ll_0,$   the PDE for $u: [0,T]\to\R^d$:
\beq\label{1POE2} (\pp_t+ L_t)u_t=\ll u_t -b_t^{(0)},\ \ t\in [0,T], u_T=0\end{equation}
has a unique solution   such that
\beq\label{AC} f_0:= \|\nn^2 u\| +|(\pp_t+\nn_{b^{(1)}})u| \in \tt L_{q_0}^{p_0},\ \  \|u\|_\infty + \|\nn u\|_\infty \le\ff 1 2.\end{equation}
 Let $\Theta_t:= id +u_t$ and
 $$\tt b_t:=\{ \ll u_t + b_t^{(1)}\}\circ \Theta_t^{-1},\ \    \tt\si_t:= \big(\{\nn\Theta_t\}\si_t\big)\circ \Theta_t^{-1},\ \ t\in [0,T].$$
 By \eqref{1POE2},     and It\^o's formula \cite[Lemma 3.3]{YZ0}, $Y_t^{x}:= \Theta_t(X_t^x)$ solves
\beq\label{YX} \d Y_t^{x} = \tt b_t(Y_t^{x})\d t+ \tt\si_t(Y_t^{x})\d W_t,\ \ Y_0^{x}= \Theta_0(x).\end{equation}
Moreover, by  $(H)$,  we find a constant $\kk>0$ such that
 \beq\label{DN10}\|\nn \tt b\|_\infty+ \|\tt\si\|_\infty+ \|(\tt\si\tt\si^*)^{-1}\|_\infty<\infty, \ \   \|\nn \tt\si\| \le \kk \sum_{i=0}^l f_i.\end{equation}
 By   \cite[Theorem 3.1]{YZ0}, for any   $(p,q)\in\scr K$, there exists a constant $c>0$ such that
$$\E\bigg(\int_s^tf_r(X_r)^2\d r\bigg|\F_s\bigg)\le c \|f\|_{\tt L_{q}^{p}(s,t)}^2,\ \ f\in \tt L_{q}^{p}(s,t),\ \ 0\le s\le t\le T,$$
where $\tt L_{q}^{p}(s,t)$ is defined as $\tt L_{q}^{p}$ for $[s,t]$ replacing $[0,T].$ This implies
Khasminskii's estimate (see for instance \cite[Lemma 4.1]{XXZZ}):  there exists an increasing  map $\Psi: (0,\infty)\times [0,\infty) \to (0,\infty)$ such that for any solution $X_t$ of \eqref{E10'},
\beq\label{KHA} \E\big[\e^{\ll \int_0^T|f_t|^2(X_t)\d t}\big]\le \Psi(\ll, \|f\|_{\tt L_q^p}),\ \ \ll>0, f\in \tt L_{q}^p,\ \ (p,q)\in \scr K.\end{equation}
Combining this with \eqref{DN1} and $f_i\in \tt L_{q_i}^{p_i}(T)$ for $(p_i,q_i)\in \scr K$, we obtain
$$\E\big[\e^{N \int_0^T \big(\|\nn \tt b\|_\infty^2 + \|\nn\tt \si_t\|^2\big)\d t}\big]<\infty,\ \ N\ge 1.$$
So,    for any $v,x\in \R^d$,   the linear SDE
\beq\label{VS0} \d v_t=  (\nn_{v_t} \tt b_t)(Y_t^{x}) +  (\nn_{v_t} \tt\si_t)(Y_t^{x}) \d W_t,\ \ v_0= v+ \nn_v u_0(x)\end{equation}
has a unique solution, and by It\^o's formula and the stochastic Gronwall inequality Lemma 3.7 in \cite{XZ},    for any $j\ge 1$ there exists a constant $c(j)>0$ such that
\beq\label{TX1} \sup_{x\in\R^d} \E \Big[\sup_{t\in [0,T]} |v_t|^j\Big] \le c(j) |v|^j,\ \ j\ge 1.\end{equation}

(b) Proof of assertion (1).  Let $Y_t^{x+\vv v}:=\Theta_t(X_t^{x+\vv v}).$
By \eqref{TX1},   for the first assertion it suffices to prove
\beq\label{DRR} \lim_{\vv\to 0} \E\Big[\sup_{t\in [0,T]} \Big|\ff{Y_t^{x+\vv v}-Y_t^{x}}\vv -v_t^{x} \Big|^j\Big] =0,\ \ j\ge 1.\end{equation}
Indeed, by an approximation argument indicated in Remark 2.1 below, see also Remark 2.1 in \cite{YZ0}, we may assume that $\nn^2 b_t^{(1)}$ is bounded so that by Lemma 2.3(3) in \cite{YZ0},
\beq\label{AC''} |\nn \Theta_t(x)-\nn \Theta_t(y)|\le c |x-y|^\aa,\ \ t\in [0,T], x,y\in \R^d\end{equation}
holds for some constants $c>0$ and $\aa\in (0,1)$. Combining this with
  \eqref{AC} and \eqref{DRR}, we see that  $\nn_v X_t^x$ exists in $L^j(\OO\to C([0,T];\R^d),\P)$ with
$$\nn_v X_t^x= (\nn\Theta_t(X_t^x))^{-1} \nn_v Y_t^x= (\nn\Theta_t(X_t^x))^{-1}v_t^x,\ \ t\in [0,T].$$
 To prove \eqref{DRR}, let
$$v_s^\vv:= \ff{Y_s^{x+\vv v}-Y_s^x}\vv,\ \ s\in [0,T], \vv \in (0,1].$$ By \eqref{DN10}, \eqref{KHA},  \cite[Lemma 2.1]{XXZZ}, and the stochastic
 Gronwall inequality   \cite[Lemma 3.7]{XZ},
 as in the proof of \cite[(4.8)]{XXZZ} we have
\beq\label{TX60}\sup_{\vv\in (0,1]}  \E\Big[\sup_{t\in [0,T]} |\tt v_t^\vv|^j\Big] <\infty,\ \ j\ge 1.\end{equation}
 Write
\beq\label{DN00}   v_r^\vv=   \int_0^r (\nn_{v_s^\vv} \tt b_s)(Y_s^x)  \d s
 +\int_0^r  (\nn_{v_s^\vv}\tt \si_s)(Y_s^x) \d W_s+\aa_r^\vv,\ \ r\in [0,t], \end{equation}
where
\beq\label{ETT}  \aa_r^\vv:=  \int_0^r \xi_s^\vv \d s+\int_0^t  \eta_s^\vv\,\d W_s  \end{equation} for
\beg{align*} &\xi_s^\vv:= \ff{\tt b_s(Y_s^{x+\vv v})- \tt b_s(Y_s^x)}\vv- (\nn_{v_s^\vv} \tt b_s)(Y_s^x),\\
&\eta_s^\vv:=    \ff{\tt\si_s(Y_s^{x+\vv v})- \tt \si_s(Y_s^x)}\vv- (\nn_{v_s^\vv} \tt \si_s)(Y_s^x).\end{align*}
We aim to prove
\beq\label{DN3}\lim_{\vv\to 0}  \E\bigg[\sup_{t\in [0,T]} |\aa_t^\vv|^n \bigg]=0,\ \ n\ge 1.\end{equation}

Firstly, since $\nn \tt b_s$ and $\nn\tt\si_s$ exist a.e., for a.e. $x\in \R^d$ we have
$$\lim_{\vv\downarrow 0} \sup_{|v|\le 1} \bigg\{ \Big|\ff{\tt b_s(x+\vv v)- \tt b_s(x)}\vv -\nn_v \tt b_s(x)\Big|+\Big\|\ff{\tt \si_s(x+\vv v)- \tt \si_s(x)}\vv -\nn_v \tt \si_s(x)\Big\|\bigg\}=0.$$
Combining this with  \eqref{TX60} and noting that
 $\L_{Y_s^x} (s\in (0,T])$ is absolutely continuous with respect to the Lebesgue measure, see for instance  Theorem 6.3.1   in \cite{BKRS},
  we obtain
\beq\label{000} \lim_{\vv\to 0} \big\{|\xi_s^\vv|+ \|\eta_s^\vv\|\big\}=0,\ \ \P\text{-a.s.},\ s\in (0,T].\end{equation}
Next,  we introduce the Hardy-Littlewood maximal functional
$$\scr M f(x):= \sup_{r\in (0,1]} \ff 1 {|B(0,r)|} \int_{B(0,r)} f(x+z)\d z,\ \ x\in \R^d, 0\le f \in L_{loc}^1(\R^d).$$
Let  $\theta>1$ such that $(\theta^{-1} p_i,\theta^{-1} q_i)\in \scr K, 0\le i\le l$. By $f_i\in \tt L_{q_i}^{p_i}$,   Lemma 2.1 in \cite{XXZZ} and  \eqref{KHA}  with $f=f_i^\theta$ and $(p,q)= (\theta^{-1} p_i, \theta^{-1}q_i)$, we obtain
\beq\label{GGR} \sup_{\vv\in [0,1]} \E\bigg[ \bigg(\int_0^T (\scr M f_i^{2\theta})(X_t^{x+\vv v}) \d t\bigg|\F_0\bigg)^n\bigg] \le K_n,\ \ 0\le i\le l\end{equation}
for some constant $K_n>0$.
By $(H)$ and Lemma 2.1 in \cite{XXZZ}, there exists a constant $c_1>0$ such that
$$|\xi_s^\vv|^{2\theta}+ \|\eta_s^\vv\|^{2\theta} \le c_1  |\tt v_t^\vv|^{2} \Big(1+\sum_{i=0}^l \big\{(\scr M f_i^{2\theta}(s, \cdot))(X_s^x) +(\scr M f_i^{2\theta}(s, \cdot))(X_s^{x+ \vv v})\big\}\Big).$$
Combining this with      \eqref{TX60} and \eqref{GGR},  for any $n\ge 1$ we find   constants $c_1(n), c_2(n)  >0$ such that
\beg{align*} &\E \bigg[\bigg(\int_0^T \big\{|\xi_s^\vv|^{2\theta}+ \|\eta_s^\vv\|^{2\theta}\big\}\d s\bigg)^n\bigg] \\
&\le c_1(n) \E\bigg[\Big(\sup_{s\in [0,T]} |v_s^\vv|^{2n} \Big) \bigg(\int_0^T \Big\{1+ \sum_{i=0}^l  (\scr M f_i^{2\theta}(s,\cdot))(X_s^x)
+(\scr M f_i^{2\theta} (s,\cdot))(X_s^{x+ \vv v})\Big\}\d s\bigg)^n\bigg]\\
&\le c_1(n) \bigg( \E\bigg[\Big(\sup_{s\in [0,T]} |v_s^\vv|^{4n} \Big) \bigg]\bigg)^{\ff 1 2}\\
&\qquad \times \bigg(\E\bigg[\bigg(\int_0^T \Big\{1+ \sum_{i=0}^l (\scr M f_i^{2\theta}(s,\cdot))(X_s^x) +(\scr M f_i^{2\theta} (s,\cdot))(X_s^{x+ \vv v})\big\}\d s\bigg)^{2n}\bigg]\bigg)^{\ff 1 2}\\
&\le c_2(n)<\infty,\ \ \vv\in (0,1].\end{align*}
Thus, by \eqref{000} and the dominated convergence theorem, we derive
 $$\lim_{\vv\to 0} \E \bigg[\bigg(\int_0^T \big\{|\xi_s^\vv|^{2}+ \|\eta_s^\vv\|^{2}\big\}\d s\bigg)^n\bigg]=0,\ \ n\ge 1.$$
Therefore,   \eqref{ETT} and      BDG's inequality imply  \eqref{DN3}.

 Finally,  by  \eqref{DN10},  \eqref{VS0},  \eqref{DN00},   and Lemma 2.1 in \cite{XXZZ},    for any $j\ge 1$, we find a  constant $c(j) >0$ such that
\beg{align*} \d |v_s-v_s^\vv|^{2j} \le  c(j)  \Big\{1+\sum_{i=0}^l f_i^2(s, Y_s^x) \Big\}|v_s-v_s^\vv|^{2j}\d s
 + c(j)  \sup_{r\in [0,s]}|\aa_r^\vv|^{2j} + \d M_s,\ \ s\in [0,t]\end{align*}
holds for some local martingale $M_s$.
Since $\lim_{\vv\to 0} |v_0 - v_0^\vv|=0$, by combining this with    \eqref{KHA}, \eqref{DN3},     and  the stochastic Gronwall inequality \cite[Lemma 3.7]{XZ}, we prove \eqref{DRR}.

 (c)  Proof of \eqref{BSM} for $f\in C_{Lip}(\R^d),$ the space of Lipschitz  continuous functions on $\R^d$. Let $t\in (0,T]$ be fixed, and consider
\beq\label{HD} h_s=\int_0^s \bb_r' \big[\tt\si_r^*\{\tt\si_r\tt\si_r^*\}^{-1} \big](Y_r^{x})v_r \d r,\ \ s\in [0,t].\end{equation}
By the same reason leading to \eqref{TX1}, the SDE
\beq\label{TX22} \beg{split}&\d w_s=\big\{ \nn_{w_s}\tt b_s(Y_s^{x}) +\tt \si_s(Y_s^{x}) h_s'\big\}\d s + (\nn_{w_s} \tt\si_s)(Y_s^{x}) \d W_s,\\
& w_0=0, s\in [0,t]\end{split} \end{equation}
has a unique solution satisfying
\beq\label{TX33} \sup_{x\in\R^d} \E \Big[\sup_{t\in [0,T]} |w_s|^j\Big] <\infty,\ \ j\ge 1.\end{equation}
We aim to prove that the Malliavin derivative $D_h Y_t^{x}$ of $Y_t^{x}$ along $h$ exists and
\beq\label{TX44} D_h Y_t^{x}= w_t.\end{equation}
For any $\vv>0$,  according to the proof of Lemma 3.1 in \cite{W21c}, \eqref{DN10},  \eqref{KHA}  and \eqref{HD} imply the well-posedness of the SDE
\beq\label{TX55} \beg{split}&\d Y_s^{x,\vv}= \big\{\tt b_t(Y_s^{x,\vv}) +\vv\tt \si_s(Y_s^{x,\vv}) h_s'\big\} \d s+ \tt\si_s(Y_s^{x,\vv}) \d W_s,\\
& s\in [0,t], Y_0^{x,\vv}=Y_0^x.\end{split}\end{equation}
By \eqref{DN10}, \eqref{HD}, Lemma 2.1 in \cite{XXZZ} and It\^o's formula, for any $j\ge 1$ we find a constant $c_1(j)>0$ such that
\beg{align*} \d |Y_s^{x,\vv}- Y_s^x|^{2j}\le &\, c_1(j)|Y_s^{x,\vv}-Y_s^x|^{2j}\sum_{i=0}^l \big\{1+\scr M \{f_i(s,\cdot)\}^2 (Y_s^x)+ \{\scr Mf_i(s,\cdot)\}^2 (Y_s^{x,\vv})\big\}\d s\\
&+ c_1(j) \vv^{2j} |v_s|^{2j} \d s +\d M_s,\ \ s\in [0,t]\end{align*}
holds for some local martingale $M_s$. Noting that $Y_0^{x,\vv}-Y_0^x=0$, by combining this with the stochastic Gronwall inequality Lemma 3.7 in \cite{XZ}  and Lemma 2.1 in \cite{XXZZ},
we obtain
\beq\label{TX66} \sup_{\vv\in (0,1]} \E \Big[\sup_{t\in [0,T]} \ff{|Y_s^{x,\vv}-Y_s^x|^j}{\vv^j}  \Big] <\infty,\ \ j\ge 1.\end{equation}
 Let $w_s^\vv= \ff{Y_s^{x,\vv}-Y_s^x}\vv.$ Then
\beq\label{DN00} \beg{split} w_r^\vv= & \int_0^r\big\{(\nn_{w_s^\vv} \tt b_s)(Y_s^x)+   \tt\si_s(Y_s^{x}) h_s'\big\} \d s\\
&+\int_0^r  (\nn_{w_s^\vv} \tt\si_s)(Y_s^x) \d W_s+\tt\aa_r^\vv,\ \ r\in [0,t]\end{split}\end{equation}
holds for
\beg{align*} \tt\aa_r^\vv:=&\ \int_0^r \Big\{\ff{\tt b_s(Y_s^{x,\vv})- \tt b_s(Y_s^x)}\vv- (\nn_{w_s^\vv} \tt b_s)(Y_s^x)\Big\}\d s\\
&+ \int_0^r \big\{\si_s(Y_s^{x,\vv})- \tt\si_s(Y_s^x)\big\} h_s'\Big\}\d s\\
&+  \int_0^r \Big\{\ff{\tt\si_s(Y_s^{x,\vv})- \tt\si_s(Y_s^x)}\vv- (\nn_{w_s^\vv} \tt\si_s)(Y_s^x)\Big\}\d W_s.\end{align*}
Combining this with \eqref{TX22} and using the same argument leading to  \eqref{DRR}, we prove   \eqref{TX44}.

  By    \eqref{HD}   and the SDE \eqref{VS0}  for $ v_s$, we see that $\bb_s v_s$ solves \eqref{TX22}, so that by the uniqueness and $\bb_t=1$, we obtain
 $$\nn_{v} Y_t^{x}= \bb_t v_t= w_t= D_h Y_t^{x}. $$
 For  $f\in C_{Lip}(\R^d)$,  $\nn f$ exists a.e. and $\|\nn f\|_\infty<\infty$. Since   $\L_{X_t^x}$ is absolutely continuous,  see for instance
   Theorem 6.3.1   in \cite{BKRS},  we conclude that $(\nn f)(X_t^x)$ exists and is bounded.
By the integration by parts formula in Malliavin calculus, see for instance   \cite{Fang}, $\nn_v Y_t^x= D_h Y_t^x$ implies
\beg{align*} &\nn_v P_t f(x)= \nn_v \E[\{f\circ (\Theta_t)^{-1}\}(Y_t^{x}) ] = \E\big[\<\nn (f\circ \Theta_t^{-1})(Y_t^{x}), \nn_v Y_t^{x}\>\big]\\
&= \E\big[D_h\{(f\circ \Theta_t^{-1})(Y_t^{x}) \}\big] = \E\bigg[f(X_t^x)\int_0^t  \< h_s',\d W_s\>\bigg]\\
&=  \E\bigg[f(X_t^x) \int_0^t \bb_s'\big\<\big\{\tt\si_s^*(\tt\si_s\tt\si_s^*)^{-1} \big\}(Y_s^{x}) v_s,  \d W_s\big\> \bigg],\ \ f\in L_{Lip}(\R^d).\end{align*}
 By $v_t=\nn Y_t^{x},$    $Y_t^{x}=\Theta_t(X_t^{x})$ and $\tt\si_t= \{(\nn\Theta_t)\si_t\}\circ \Theta_t^{-1}$,  we obtain
\beg{align*} &\{\tt\si_s^*(\tt\si_s\tt\si_s^*)^{-1} \big\}(Y_s^{x}) v_s= \Big[\si_s^*(\si_s\si_s^*)^{-1}\big\{(\nn\Theta_s) \si_s\si_s^*(\nn\Theta_s)^*\big\}^{-1}\Big](X_s^x) \big\{\nn\Theta_s(X_s^x)\big\}\nn X_s^x\\
 &= \{\si_s^*(\si_s\si_s^*)^{-1} \big\}(X_s^x) \nn X_s^x,\ \ s\in [0,T],\end{align*}
  so that this implies
\beq\label{BSM'}  \nn_v P_tf(x)= \E\bigg[f(X_t^x) \int_0^t \bb_s'\big\<\big\{\si_s^*(\si_s\si_s^*)^{-1} \big\}(X_s^x) \nn X_s^x, \d W_s\big\> \bigg],\ \ f\in C_{Lip}(\R^d).\end{equation}

(d) Proof of \eqref{GRD} and \eqref{BSM}. Let $P_t^*\dd_x=\L_{X_t^x}$ and let $\nu_\vv$ be the finite signed measure defined by
$$\nu_\vv(A):=  \int_0^\vv  \E\bigg[1_A(X_t^{x+r v}) \int_0^t \bb_s'\big\<\big\{\si_s^*(\si_s(\si_s)^*)^{-1} \big\}(X_s^{x+rv}) \nn X_s^{x+rv}, \d W_s\big\> \bigg]\d r$$
for $A\in \B(\R^d),$ the Borel $\sigma$-algebra on $\R^d$.
Then \eqref{BSM'} implies
$$ (P_t^* \dd_{x+\vv v} -P_t^*\dd_x)(f)= \nu_\vv(f),\ \ f\in C_{Lip}(\R^d),$$
where $\nu(f):=\int f\d \nu$ for a (signed)  measure $\nu$ and $f\in L^1(|\nu|)$.
Since $C_{Lip}(\R^d)$ determines measures, we obtain
$$P_t^* \dd_{x+\vv v} -P_t^*\dd_x= \nu_\vv,$$
so that for any $f\in \B_b(\R^d)$,
\beg{align*} &P_tf(x+\vv v)- P_t f(x)\\
&= \int_0^\vv   \E\bigg[f(X_t^{x+r v}) \int_0^t \bb_s'\big\<\big\{\si_s^*(\si_s(\si_s)^*)^{-1} \big\}(X_s^{x+rv}) \nn X_s^{x+rv}x, \d W_s\big\> \bigg]\d r.\end{align*}
Combining this with \eqref{DFF0} and the boundedness of $\si^*(\si\si^*)^{-1}$, we prove \eqref{GRD}.

Next,  let $f\in \B_b(\R^d). $ For any $r\in (0,T),$ let $(X^x_{r,t})_{t\in [r,T]}$ solve \eqref{E10'} from time $r$ with $X_{r,r}^x=x$. Let
$$P_{r,t}f(x):= \E[f(X_{r,t}^x)],\ \ f\in \B_b(\R^d), x\in\R^d.$$ Then the  well-posedness implies
$$P_t= P_rP_{r,t},\ \ 0<r<t\le T.$$
Moreover, considering the SDE from time $r$ replacing $0$, \eqref{GRD} implies
$$\|\nn P_{r,t}f\|_\infty<\infty,\ \ f\in \B_b(\R^d), 0<r<t\le T.$$
So,  by \eqref{BSM'}  for $(P_r,\bb_s/\bb_r)$ replacing $(P_t,\bb_s)$, we obtain
\beg{align*} & \nn_v P_tf(x)=\nn_vP_r(P_{r,t}f)(x)\\
&= \ff 1 {\bb_r} \E\bigg[P_{r,t}f(X_r^x) \int_0^r \bb_s'\big\<\big\{\si_s^*(\si_s(\si_s)^*)^{-1} \big\}(X_s^x) \nn X_s^x, \d W_s\big\> \bigg] \end{align*}
for all $f\in\B_b(\R^d)$ and $r\in (0,t)$ such that $\bb_r>0.$  Since  the Markov property implies
$$\E[f(X_t^x)|\F_r) =P_{r,t}f(X_r^x),$$
we obtain
$$ \nn_v P_tf(x)= \ff 1 {\bb_r} \E\bigg[f(X_t^x) \int_0^r \bb_s'\big\<\big\{\si_s^*(\si_s(\si_s)^*)^{-1} \big\}(X_s^x) \nn X_s^x, \d W_s\big\> \bigg], $$
so that letting $r\uparrow t$ gives \eqref{BSM}. \end{proof}

 To conclude this section, we make the following remark which enables us to apply Theorem \ref{T03} to the decoupled SDE \eqref{E102} with estimates uniformly in $\mu$.

\paragraph{Remark 2.1.} For fixed $\si$ but may be variable $b$, the constants $c(\cdot)$ in  Theorem \ref{T03} are uniformly in  $b=b^{(0)}+b^{(1)}$ satisfying
 \beq\label{BNN} \|b^{(0)}\|_{\tt L_{q_0}^{p_0} } + \|\nn b^{(1)}\|_\infty \le N\end{equation}  for a given constant $N>0$. Indeed,
letting $\gg$ be the standard Gaussian measure and take
$$\tt b_t^{(1)} (x):= \int_{\R^d} b_t^{(1)}(x+y) \gg(\d y),\ \ x\in \R^d, t\in [0,T],$$
we find   constant $c >0$ only depending on $N$  such that \eqref{BNN} implies
$$\|\nn \tt b^{(1)}\|_\infty + \|\nn^2 \tt b^{(1)}\|_\infty + \|b_t^{(1)}- \tt b_t^{(1)}\|_\infty\le c.$$
Then $\tt b^{(0)}:= b^{(0)}+ \tt b^{(1)}- b^{(1)}$ satisfies
$$\|\tt b^{(0)}\|_{\tt L_{q_0}^{p_0}}\le \|b^{(0)}\|_{\tt L_{q_0}^{p_0}} + c \| 1\|_{\tt L_{q_0}^{p_0}}=:c'.$$
 According to the proofs of \cite[Theorem 2.1 and Theorem 3.1]{YZ0}  for $b= \tt b^{(0)} +\tt b^{(1)}$,  the  constant  $\ll_0>0$ before  \eqref{1POE2}, $\|\nn^2 u\|_{\tt L_{q_0}^{p_0}}$  and the constant in Krylov's
 estimate are uniformly in $b$ satisfying \eqref{BNN}.
 According to the proof of  \cite[Lemma 4.1]{XXZZ}, the same is true for  Khasminskii's estimate \eqref{KHA}.
Therefore, in the proof of Theorem \ref{T03},  constants $c(\cdot)$ can be taken uniformly in $b$  satisfying \eqref{BNN}.

\section{Proof of Theorem \ref{TA2} }

We first present a Lipschitz estimate for $L$-differentiable functions on $\scr P_k$.

 \beg{lem}\label{LL2} Let  $f$ be $L$-differentiable on $\scr P_k$   such that for any $\mu\in \scr P_k$,  $D^L f(\mu)(\cdot)$ has a continuous version satisfying
\beq\label{PLL} |D^Lf(\mu)(x)|\le c(\mu) (1+|x|^{k-1}),\ \ x\in\R^d \end{equation}
 holds for come constant  $c(\mu)>0$, and
\beq\label{FF1} K_0:=\sup_{\mu\in \scr P_k} \|D^L f(\mu)\|_{L^{k^*}(\mu)} <\infty.\end{equation}  Then
\beq\label{FFE} |f(\mu_1)-f(\mu_2)|\le K_0 \W_k(\mu_1,\mu_2) ,\ \ \mu_1,\mu_2\in \scr P_k.\end{equation}
\end{lem}

\beg{proof} Let $\xi_1,\xi_2$ be two random variables with
$$\L_{\xi_1}=\mu_1,\ \ \L_{\xi_2}=\mu_2,\ \ \W_k(\mu_1,\mu_2)=(\E[|\xi_1-\xi_2|^k])^{\ff 1 k}.$$
Let $\eta$ be a normal random variable independent of $(\xi_1,\xi_2)$. Then
$$\gg_{\vv}(r):= \vv \eta + r\xi_1+ (1-r)\xi_2,\ \ r\in [0,1], \vv\in (0,1]$$ are absolutely continuous with respect to the Lebesgue measure and hence atomless.
 By Theorem 2.1 in \cite{BRW20}, \eqref{PLL} and the continuity of $D^Lf(\mu)(\cdot)$  imply
\beg{align*} &|f(\L_{\gg_\vv(1)})- f(\L_{\gg_\vv(0)})| = \bigg|\int_0^1 \E\big[\<D^L f(\L_{\gg_\vv(r)})(\gg_\vv(r)), \xi_1-\xi_2\>\big]\d r\bigg|\\
&\le \big(\E[|\xi_1-\xi_2|^k]\big)^{\ff 1 k}\int_0^1 \|D^L f(\L_{\gg_\vv(r)})\|_{L^{k^*}( \L_{\gg_\vv(r)})}\d r\\
&\le K \W_k(\mu_1,\mu_2),\ \ \vv\in (0,1].\end{align*}
Letting $\vv\to 0$ we prove \eqref{FFE}.
\end{proof}

 In the following, we prove assertions (1) and (2) by using Zvonkin's transform.

\beg{proof}[Proof of Theorem \ref{TA2}(1)]    By $(H)$, we have $b^{(1)}_t \in \D_k$ with $\|D^L b^{(1)}_t(x,\mu)\|_{L^{k^*}(\mu)} \le K$ for some constant $K>0.$ Then    Lemma \ref{LL2}   implies
\beq\label{BKK} |b^{(1)}_t(x,\mu)- b^{(1)}_t(x,\nu)|\le K \W_k(\mu,\nu),  \end{equation}
  so that  the well-posedness of $\eqref{E10}$ follows from
Theorem 3.1(2)(ii) in \cite{W21b} for $D=\R^d$ for which condition (3) can be dropped from $(A_2^{\si,b})$ therein since $\pp D=\emptyset$, it is also implied by Theorem 1.1(2) in \cite{HW20} where in   (9) the condition  $\bar b^\mu\in L_{q}^p$ can be weakened as $\bar b^\mu\in \tt L_{q}^p$, since in the proof we may replace $L_{q}^p$ by $\tt L_q^p$ according to   Theorem 2.1 and Theorem 3.1  in \cite{YZ0}.

To prove \eqref{EST0} and \eqref{*BN}, we use Zvonkin's transform.
 Consider the differential operator
\beq\label{PDEN} L_t^\mu=\ff 1 2{\rm tr}\big\{\si_t\si_t^*\nn^2\big\}+\nn_{b_t(\cdot,\mu_t)},\ \ t\in [0,T].\end{equation}
By \cite[Theorem 2.1]{YZ0}, $(H)$ implies that for some $\ll_0$ uniformly in $\mu_0$, when $\ll\ge \ll_0$  the PDE
\beq\label{PDEM} (\pp_t+ L_t^\mu)u_t=\ll u_t -b_t^{(0)},\ \ t\in [0,T], u_T=0\end{equation}
has a unique solution $u\in \tt H_{q_0}^{2,p_0}$ such that \eqref{AC} holds.
 Let $\Theta_t:= id +u_t.$ By It\^o's formula in Lemma 3.3 of \cite{YZ0},
 $$Y_t:= \Theta_t(X_t)= X_t+u_t(X_t)$$
 solves the SDE
 \beq\label{*Y0}    \d Y_t =    \big\{b^{(1)}_t(X_t,\mu_t)+\ll u_t(X_t)\big\}    \d t
 + \big\{(\nn\Theta_t) \si_t\big\} (X_t)   \d W_t,\ \ Y_0= \Theta_0(X_0).
  \end{equation}
By  \eqref{AC}, there exists a constant  $c_1>1$ such that
\beq\label{CCA} |X_t| \le  c_1 (1+|Y_t|)\le c_1^2 (1+ |X_t|),\ \ t\in [0,T].\end{equation}
For any $n\ge 1$, let $$\gg_{t,n}:= \sup_{s\in [0,t\land\tau_n]} |Y_s|,\ \ \tau_n:=\inf\{s\ge 0: |Y_s|\ge n\}, t\in [0,T].$$
By BDG's inequality,  $(H)$ and \eqref{AC}, for any $j\ge 1$ there exists a constant $c(j)>0$ such that
$$  \E\big(\gg_{t,n}^j\big|\F_0\big) \le 2 |Y_{0}|^j + c(j) \int_0^t \big\{ \E(\gg_{s,n}^j|\F_0)+ (\E[|Y_s|^k])^{\ff j k} +1 \big\}\d s + c(j),\ \ n\ge 1, t\in [0,T].$$
By Gronwall's inequality,
\beq\label{AA01} \E\big(\gg_{t,n}^j\big|\F_0\big) \le \bigg(2 |Y_{0}|^j + c(j) \int_0^t \big\{  (\E[|Y_s|^k])^{\ff j k}   +1 \big\}\d s + c(j)\bigg)\e^{c(j)t},\ \ n\ge 1, t\in [0,T]. \end{equation}
Taking expectations with $j=k$ and letting $n\to\infty$, we find a constant $c_2>0$ such that
$$\E[\gg_t^k]:= \E\bigg[\sup_{s\in [0,t]} |Y_s|^k\bigg]\le c_2 (1+ \E[|Y_0|^k]) + c_2\int_0^t \E[|Y_s|^k]\d s,\ \ t\in [0,T].$$
Noting that $\sup_{t\in [0,T]} \E[|X_t|^k]<\infty$ as $X_t$ is the solution of \eqref{E10} for distributions in $\scr P_k$, by combining this with \eqref{CCA} and
$\E[\gg_t^k]\ge \E[|Y_s|^k]$ we obtain
$$\E[\gg_t^k]:= \E\bigg[\sup_{s\in [0,t]} |Y_s|^k\bigg]\le c_2 + c_2\int_0^t \E[\gg_s^k]\d s<\infty,\ \ t\in [0,T],$$
so that by Gronwall's inequality and \eqref{CCA}, we prove  \eqref{*BN} for some constant $c>0.$  Substituting this into \eqref{AA01} and letting $n\to\infty$, we
prove \eqref{EST0}.  \end{proof}

\beg{proof}[Proof  of Theorem \ref{TA2}(2)]  Denote $\mu_t^i:=\L_{X_t^i}, i=1,2, t\in [0,T]$. Let $u$ solve \eqref{PDEM} for $L_t^{\mu^1}$ replacing $L^\mu_t$ such that \eqref{AC} holds.
Let $\Theta_t= id +u_t$ and
$$Y_t^i= \Theta_t(X_t^i),\ \ t\in [0,T], i=1,2.$$
By \eqref{PDEM} and It\^o's formula we obtain
\beg{align*} &\d Y_t^1 = \big\{b_t^{(1)}(X_t^1, \mu_t^1) +\ll u_t(X_t^1)\big\}\d t +\big\{(\nn\Theta_t)\si_t\big\}(X_t^1) \d W_t,\\
&\d Y_t^2 = \big\{b_t^{(1)}(X_t^2, \mu_t^2) +\ll u_t(X_t^2)+\nn_{b_t^{(1)}(X_t^2,\mu_t^2)- b_t^{(1)}(X_t^2,\mu_t^1) } u_t(X_t^2)\big\}\d t \\
&\qquad \qquad +\big\{(\nn\Theta_t)\si_t\big\}(X_t^2) \d W_t,\ \ t\in [0,T].\end{align*}
So, by It\^o's formula, the process
$$v_t:=  Y_t^2-Y_t^1, \ \ \ t\in [0,T]$$ satisfies the SDE
  \beg{align*}   \d v_t=  &\,   \Big\{b^{(1)}_t(X_t^2, \mu_t^2)+\ll u_t(X_t^2)- b_t^{(1)}(X_t^1,\mu_t^1) -\ll u_t(X_t^1)  +\nn_{b^{(1)}_t(X_t^2,\mu_t^2)- b^{(1)}_t (X_t^2, \mu_t^1)}u_t(X_t^2)\Big\}    \d t \\
 &+ \Big\{(\nn\Theta_t) \si_t\big\} (X_t^2) - (\nn\Theta_t) \si_t\big\} (X_t^1)\Big\}  \d W_t,\ \ v_0 = \Theta_0(X_0^2)-\Theta_0(X_0^1).
 \end{align*}
 By \eqref{AC} and \eqref{BKK}, we obtain
 $$|b_t^{(1)}(x,\mu_t^2)- b_t^{(1)}(x,\mu_t^1)|^k\le K^k \E[|X_t^2- X_t^1|^k]\le (2K)^k \E[|Y_t^2-Y_t^1|^k].$$
Combining this with  $(H)$, \eqref{AC},      Lemma 2.1 in \cite{XXZZ},  and applying It\^o's formula,  for any $j\ge k$ we find a constant $c_1>0$ such that
\beq\label{TTYU}  \beg{split} |v_t|^{2j} \le &\, |v_0|^{2j} + c_1 \int_0^t |v_s|^{2j}\Big\{1+\sum_{i=0}^l \scr M f_i^2(s,X_s)\Big\} \d s\\
&+c_1\int_0^t  (\E[|v_s|^k])^{\ff{2 j}k}  \d s+ M_t,\ \ t\in [0,T]\end{split} \end{equation}
holds for some local martingale $M_t$ with $M_0=0.$  Since \eqref{AC} implies
$$|v_0|\le   2 |X_0^1-X_0^2|,$$
by stochastic Gronwall's inequality  \cite[Lemma 3.7]{XZ},  Lemma 2.1 in \cite{XXZZ} and Khasminskii's estimate \eqref{KHA}, we find a constant $c_2>0$ such that
$$\gg_t := \sup_{s\in [0,t]} |v_s|,\ \ t\in [0,T]$$ satisfies
\beq\label{LOP} \beg{split} & \E\big[|\gg_t|^j \big|\F_0\big] \le c_2 \bigg( |X_0^1-X_0^2|^{2j} + \int_0^t  (\E[|v_s|^k])^{\ff{2j}k} \d s\bigg)^{\ff 1 2} \\
&\le c_2|X_0^1-X_0^2|^j + \ff 1 2 \sup_{s\in [0,t]} \big(\E[|v_s|^k]\big)^{\ff j k} + \ff {c_2^2}2 \int_0^t  \big(\E[|v_s|^k]\big)^{\ff j k} \d s <\infty,\ \ t\in [0,T].\end{split}\end{equation}
Noting that $\sup_{s\in [0,t]} \E[|v_s|^k] \le \E[|\gg_t|^k]$, by taking    expectation  in \eqref{LOP} with $j=k$,  we derive
$$\E\big[|\gg_t|^k\big]\le 2 c_2 \E[|X_0^1-X_0^2|^k] + c_2^2 \int_0^t \E[|\gg_s|^k]\d s,\ \ t\in [0,T].$$
Since $ \E\big[|\gg_t|^k\big]<\infty$ due to \eqref{*BN}, by Gronwall's inequality we find a constant $c>0$ such that
$$\sup_{t\in [0,T]} \E[|v_s|^k]\le \E[|\gg_T|^k]\le c \E[|X_0^1-X_0^2|^k].$$
  Substituting this into \eqref{LOP}  implies \eqref{EE0}.  \end{proof}

\beg{proof}[Proof  of Theorem \ref{TA2}(3)]   Let $\nu\in \scr P_k$ and take $\F_0$-measurable random variables $X_0, \tt X_0$ such that
\beq\label{N1} \L_{X_0}=\mu, \ \ \L_{\tt X_0} =\nu,\ \ \E[|X_0-\tt X_0|^k]= \W_k(\mu,\nu)^k.\end{equation}
Let $X_t$ and $\tt X_t$ solve \eqref{E10} with initial values $X_0$ and $ \tt X_0$ respectively, and denote
$$ \mu_t:= P_t^*\mu= \L_{X_t},\ \ \nu_t:= P_t^*\nu=\L_{\tt X_t},\ \ t\in [0,T]. $$
  Let $P_t^\mu$ be the semigroup associated with $X_t^{\mu,x}$.  According to Remark 2.1, \eqref{GRD} holds for  $P_t^\mu$ replacing $P_t$ and  some  constant $c>0$ independent of $\mu$. Then
  \beq\label{LN3} \|P_t^*\mu - (P_t^\mu)^*\nu\|_{var}= \|(P_t^\mu)^*\mu - (P_t^\mu)^*\nu\|_{var} \le \ff{c}{\ss t} \W_1(\mu,\nu).\end{equation}
  On the other hand,  let
  $$R_t:= \e^{-\int_0^t \<\zeta_s(X_s)\{b_s(\tt X_s,\nu_s)- b_s(\tt X_s,\mu_s)\}, dW_s\>- \ff 1 2 \int_0^s |\zeta_s(X_s)\{b_s(\tt X_s,\nu_s)- b_s(\tt X_s,\mu_s)\}|^2\d s}.$$
  By $(H)$ and Girsanov's theorem, $\Q_t:= R_t\P$ is a probability measure under which
  $$\tt W_s:= W_s +\int_0^s\zeta_s(X_s)\{b_s(\tt X_s,\nu_s)- b_s(\tt X_s,\mu_s)\}\d s,\ \ r\in [0,t]$$ is a Brownian  motion.  Reformulating the SDE for $\tt X_s$ as
  $$\d \tt X_s = b_s(\tt X_s, \mu_s) \d s + \si_s(\tt X_s) \d\tt W_s,\ \ \L_{\tt X_0}=\nu,$$
  by the uniqueness we obtain $\L_{\tt X_t|\Q_t}= (P_t^\mu)^*\nu,$ so that by Pinsker's inequality and $(H)$, we find constants $c_1 >0$ such that
  \beg{align*} &\|P_t^*\nu - (P_t^\mu)^{*} \nu\|_{var}^2=\sup_{|f|\le 1} |\E[f(\tt X_t)(R_t-1)]|^2\le 2 \E[ R_t\log R_t]\\
  &\le c_1 \E_{\Q_t } \int_0^t \W_k(\mu_s,\nu_s)^2\d s = c_1 \int_0^t  \W_k(\mu_s,\nu_s)^2\d s.\end{align*}
  Combining this with \eqref{LN2} and \eqref{LN3}, we prove \eqref{GRD00} for some constant $c>0$.
  \end{proof}

\section{Proof of Theorem \ref{TA2'}}

A key step of the proof is to calculate $\nn_\phi X_t^\mu$. In general, let $ X_t^\mu$ solve \eqref{E10} for $\L_{X_0^\mu}=\mu\in \scr P_k$,
and for any $\vv\in [0,1]$ and $\F_0$-measurable random variable  $\eta$ with $\L_{\eta}\in \scr P_k$,
let $X_t^\vv$ solve \eqref{E10} with $X_0^\vv= X_0^\mu+\vv\eta$.  We intend to calculate
\beq\label{ETAR} \nn_\eta X_t^\mu:= \lim_{ \vv \downarrow 0} \ff{X_t^{\vv}-X_t^\mu} \vv,\ \  t\in [0,T]\end{equation}
 in $L^k(\OO\to C([0,T]; \R^d),\P)$. In particular, taking $\eta:=\phi(X_0^\mu)$ for $\phi\in T_{\mu,k}$, we have
\beq\label{ETAR2} \nn_\phi X_t^\mu= \nn_\eta X_t^\mu,\ \ t\in [0,T].\end{equation} Choosing general $\eta$ instead of $\phi(X_0^\mu)$ is useful in the proof of Theorem \ref{TA2''}.

\subsection{ The SDE for $\nn_\eta X_t^\mu$}

Comparing with $\nn_v X_t^x$ in Section 2, there are two essential differences  in the study of $\nn_\eta X_t^\mu:$
\beg{enumerate} \item[(1)] the $L$-derivative of $b^{(1)}(x,\mu)$ will be involved;
\item[(2)]   since $X_0^\mu$ is a random variable with $\E[|X_0^\mu|^k]<\infty$,
in general we do not have $$\sup_{\vv\in [0,1]}\E[|X_t^{\vv}|^j]<\infty,\ \ j>k,$$ which is   important for the dominated convergence theorem  as used   in the study of $\nn_v X_t^x$. \end{enumerate}

Point (1)   will be managed by a chain rule
due to \cite{BRW20} for functions on $\scr P_k$, see  \cite{HSS} for earlier versions on $\scr P_2$. Point (2)    will  be treated using  the conditional expectation $\E[\cdot|\F_0]$ to replace the expectation $\E$, since
we can prove $\E[|X_t^\vv|^j|\F_0]<\infty$   for any $j\ge 1.$

Let   $u$ solve \eqref{PDEM}  such that \eqref{AC} and  \eqref{AC''} hold as explained before.   Let
$\Theta_t= id+ u_t$   and
\beq\label{YTV} Y_t^r:=\Theta_t(X_t^r)= X_t^r + u_t(X_t^r),\ \ t\in [0,T], \ r\in [0,1].\end{equation}
 By \eqref{PDEM} and It\^o's formula, see Lemma 3.3 in \cite{YZ0}, for any $r \in [0,1]$ we have
 \beq\label{*YY} \beg{split}   \d Y_t^r &=    \big\{b^{(1)}_t(X_t^r,\mu_t^r)+\ll u_t(X_t^r)+ \nn_{b^{(1)}_t(X_t^r,\mu_t^r)- b^{(1)}_t (X_t^r, \mu_t)}u_t(X_t^r)\big\}    \d t \\
 &+ \big\{(\nn\Theta_t) \si_t\big\} (X_t^r)   \d W_t,\ \ Y_0^r= \Theta_0(X_0^r)=X_0^\mu+r\eta+u_0(X_0^\mu+r\eta).
 \end{split}\end{equation}
For any $t\in [0,T]$ and $v\in L^k(\OO\to \R^d,\P),$ let
\beq\label{PSV0}  \psi_t(v):= \E\big[\<D^L   b_t^{(1)}(z, \mu_t)( X_t^\mu), (\nn\Theta_t( X_t^\mu))^{-1} v\>\big]\big|_{z= X_t^\mu}.\end{equation}
By $(H)(3)$, there exists a constant $K>0$ such that for any $ v,\tt v\in L^k(\OO\to \R^d,\P),$
\beq\label{PSV}  \psi_t(0)=0,\ \   | \psi_t(v)- \psi_t(\tt v)|\le K (\E[|v-\tt v|^k])^{\ff 1 k},\ \ t\in [0,T].\end{equation}
If
\beq\label{VLI} v_t^\eta:= \nn_\eta Y_t^0:= \lim_{\vv\downarrow 0} \ff {Y_t^\vv-Y_t^0}\vv\end{equation}  exists in $L^k(\OO\to C([0,T];\R^d),\P)$,  by \eqref{AC}, \eqref{AC''} and \eqref{YTV} we see  that $\nn_\eta  X_t^\mu$ exists in the same sense and
\beq\label{PX} \nn_\eta  X_t^\mu= (\nn\Theta_t( X_t^\mu))^{-1} \nn_\eta Y_t^0= (\nn\Theta_t( X_t^\mu))^{-1} v_t^\eta.\end{equation}
Combining this  with  $(H)$,  applying the chain rule Theorem 2.1 in \cite{BRW20},  and noting that $\mu_t$ is absolutely continuous due to  Theorem 6.3.1   in \cite{BKRS}, we obtain
\beq\label{*0*}\beg{split}   &\lim_{\vv\to 0} \ff{b_t^{(1)}(X_t^\vv,\mu_t^\vv)- b_t^{(1)}(X_t^\vv,\mu_t)}{\vv}=\lim_{\vv\to 0} \int_0^1 \ff 1 \vv \ff{\d}{\d r}  b_t^{(1)}(X_t^\vv,\L_{r X_t^\vv+ (1-r)X_t^\mu}) \d r\\
&=  \lim_{\vv\to 0} \int_0^1 \E\Big[\Big\<D^L b_t^{(1)}(z,\L_{r X_t^\vv+ (1-r)X_t^\mu} )(r X_t^\vv+ (1-r)X_t^\mu),\ff{X_t^\vv-X_t^\mu}\vv\Big\>\d r\Big]\Big|_{z=X_t^\vv}\\
&=   \psi_t(v_t^\eta),\end{split} \end{equation} which together with \eqref{PX} yields
 \beg{align*}&\lim_{\vv\to 0} \ff{b^{(1)}_t(X_t^\vv,\mu_t^\vv)- b^{(1)}_t ( X_t^\mu, \mu_t)}{\vv}=  \psi_t(v_t^\eta)+ \nn_{(\nn\Theta_t( X_t^\mu))^{-1} v_t^\eta} b_t^{(1)}( X_t^\mu,\mu_t),\\
&\lim_{\vv\to 0} \ff{ \{(\nn\Theta_t) \si_t\} (X_t^\vv)- \{(\nn\Theta_t) \si_t\} ( X_t^\mu)}\vv = \nn_{(\nn\Theta_t( X_t^\mu))^{-1} v_t^\eta}  \{(\nn\Theta_t) \si_t\} ( X_t^\mu),\\
& \lim_{\vv\to 0} \ff{u_t(X_t^\vv)- u_t ( X_t^\mu)}{\vv}=  \nn_{(\nn\Theta_t( X_t^\mu))^{-1} v_t^\eta} u_t( X_t^\mu).\end{align*}
Thus, if    $v_t^\eta$  in \eqref{VLI} exists, by \eqref{*YY} it should solve the SDE
 \beq\label{VS}\beg{split}  \d v_t^\eta= &\, \big\{\psi_t(v_t^\eta)+ \nn_{(\nn\Theta_t( X_t^\mu))^{-1} v_t^\eta} b_t^{(1)}( X_t^\mu,\mu_t) +\nn_{\psi_t(v_t^\eta)+\ll (\nn\Theta_t( X_t^\mu))^{-1}v_t^\eta} u_t( X_t^\mu)\big\} \d t   \\
&+   \nn_{(\nn\Theta_t( X_t^\mu))^{-1} v_t^\eta}  \{(\nn\Theta_t) \si_t\} ( X_t^\mu)   \d W_t,\ \ v_0^\eta= \eta+ (\nn_{\eta}  u_0)(X_0).\end{split} \end{equation}
Therefore, in terms of \eqref{PX}, to study $\nn_\eta X_t^\mu$  we first consider the SDE \eqref{VS}.

\beg{lem}\label{LL1} Assume   $(H)$.  For any   $\eta\in L^k(\OO\to\R^d,\F_0,\P)$,   the SDE
$\eqref{VS}$ has a unique solution, and for any $j\ge 1$ there exists a constant $c>0$  such that
\beq\label{EST1} \E\bigg[\sup_{t\in [0,T]} |v_t^\eta|^j\bigg|\F_0 \bigg]\le c\big\{(\E[|\eta|^k])^{\ff jk}+  |\eta|^j\big\}, \ \ \mu\in \scr P_k, \eta\in L^k(\OO\to\R^d,\F_0,\P).\end{equation}
\end{lem}

\beg{proof} We simply denote $X_t=X_t^\mu, t\in [0,T].$

 (1) Well-posedness of \eqref{VS}. Consider the space
$$\C_k:= \bigg\{(v_t)_{t\in [0,T]}\ \text{is\ continuous\ adapted,}\ v_0=v_0^\eta,\  \E\Big[\sup_{t\in [0,T]} |v_t|^k\Big]<\infty\bigg\},$$ which is complete under the metric
$$\rr_\ll(v^1,v^2):=\bigg( \E\bigg[\sup_{t\in [0,T]} \e^{-\ll t}|v_t^1-v_t^2|^k\bigg]\bigg)^{\ff 1 k},\ \ v^1,v^2\in \C_k$$ for $\ll>0$.
By $(H)$,   \eqref{AC} and \eqref{PSV}, there exist a constant $K>0$ and  a function $1\le  f_0\in\tt L_{q_0}^{p_0}$ such that for any random variable  $v$,
\beq\label{DN1} \beg{split} & \big|\nn_{(\nn\Theta_t(X_t))^{-1} v} b_t^{(1)}(X_t,\mu_t) +\nn_{\psi_t(v_t)+\ll (\nn\Theta_t(X_t))^{-1} v}  u_t(X_t) \big|\le K |v|,\\
&\big\| \nn_{(\nn\Theta_t(X_t))^{-1} v}  \{(\nn\Theta_t) \si_t\} (X_t) \big\|\le K |v| \sum_{i=0}^l f_i(t, X_t),\\
&\big|\psi_t(v)\big|\le K \big(\E[|v|^k]\big)^{\ff 1 k},\ \ t\in [0,T].\end{split}\end{equation}
Let $f=\sum_{i=0}^l f_i$. Let $\theta>1$ such that $(\theta^{-1} p_i,\theta^{-1} q_i)\in \scr K, 0\le i\le l$. By Krylov's estimate Lemma 3.2(1) in \cite{YZ0}, we find a constant $c>0$ such that
$$\E\int_0^T f_t(X_t)^{2\theta}\d t\le  c \sum_{i=0}^l \|f_i\|_{\tt L_{q_i}^{p_i}}^{2\theta}<\infty.$$
So,
$$\tau_n:=T\land \inf\bigg\{t\ge 0: \int_0^t |f_t(X_t)|^{2\theta}\d s\ge n\bigg\}\to T\ \text{as}\ n\to\infty.$$
Thus,
\beq\label{DN2}  \beg{split} H_t(v):= &\ v_0^\eta +\int_0^t  \Big\{\psi_s(v_s)+ \nn_{(\nn\Theta_s(X_s))^{-1} v_s} b_s^{(1)}(X_s,\mu_s) \\
&\qquad\qquad \qquad +\nn_{\psi_s(v_s)+\ll (\nn\Theta_s(X_s))^{-1} v_s} u_s(X_s)\Big\}   \d s\\
& + \int_0^t  \nn_{(\nn\Theta_s(X_s))^{-1} v_s}  \{(\nn\Theta_s) \si_s\} (X_s)  \d W_s,\ \ t\in [0,T]\end{split}\end{equation}
is an adapted continuous process on $\R^d$, and for any $n\ge 1$,
$$H_{\cdot\land\tau_n}: \C_{k,n}\to \C_{k,n},\ \ \C_{k,n}:= \big\{(v_{\cdot\land \tau_n}):\ \ v\in \C_k\big\}.$$
So,  it remains  to prove that $H$ has a unique fixed point $v^\eta\in \C_k$ satisfying \eqref{EST1}, which is then the unique solution of \eqref{VS}.
In the following we explain that   it suffices to prove
\beq\label{KKL} H_{\cdot\land\tau_n}\ \text{ has\ a\ unique\ fixed\ point\ in\ }    \C_{k,n},\ \   n\ge 1.\end{equation}

 Indeed, if \eqref{KKL} holds, then the unique fixed point  $v^{\eta,n}_{\cdot\land\tau_n}$ satisfies
$$v^{\eta,n}_{\cdot\land\tau_n}= v^{\eta,n+k}_{\cdot\land\tau_n},\ \ n,k\ge 1,$$
so that
$$  v_t^\eta:= \lim_{n\to\infty} v_{t\land\tau_n}^{\eta,n}$$ is a continuous adapted process on $\R^d$, and
$$ H_{\cdot\land\tau_n}( v^\eta)=  v^\eta_{\cdot\land\tau_n}\in \C_{k,n},\ \ n\ge 1.$$
By this and \eqref{DN1}, for any $j\ge k$ we find a constant $c>0$ such that
$$\d |v_t^\eta|^{2j} \le c\big(\{\E[|v^\eta_{t\land\tau_n}|^k]\}^{\ff {2j } k}  + |v^\eta_t|^{2j} \big)(1+f_t^2(X_t))\d t + \d\tt M_t,\ \ t\in [0,\tau_n],$$
holds for some local martingale $\tt M_t$. By the stochastic Gronwall inequality,  we find  constants $k_1,k_2>0$ such that
\beq\label{RS} \beg{split} \E\Big[\sup_{t\in [0,T]} |v^\eta_{t\land\tau_n}|^j\Big|\F_0\Big]
&\le k_1\bigg(\int_0^t\big\{\E[|v^\eta_{t\land\tau_n}|^k] \big\}^{\ff{2j}k} \d s+   \E[|v_0^\eta|^{2j}|\F_0]\bigg)^{\ff 1 2}\\
&\le k_2 |\eta|^{j} +k_1\bigg(\int_0^t\big\{\E[|v^\eta_{t\land\tau_n}|^k] \big\}^{\ff{2j}k} \d s\bigg)^{\ff 1 2}.\end{split}\end{equation}
Taking $j=k$ we obtain
 $$\E\Big[\sup_{t\in [0,T]} |v^\eta_{t\land\tau_n}|^k\Big|\F_0\Big] \le k_2 |\eta|^k + \ff {k_1^2} 2  \int_0^t \E[|v^\eta_{t\land\tau_n}|^k]  \d s+\ff 1 2 \E\Big[\sup_{t\in [0,T]} |v^\eta_{t\land\tau_n}|^k\Big],\ \ t\in [0,T].$$
Taking expectation and applying Gronwall's inequality, we find a constant $k_3>0$ such that
$$\sup_{n\ge 1} \E\Big[\sup_{t\in [0,T]} |v^\eta_{t\land\tau_n}|^k\Big]\le k_3 \E[|\eta|^k],$$
so that \eqref{RS} with $n\to\infty$ implies  \eqref{EST1}, and it is the unique solution of  \eqref{VS} in $\C_k$ since for each $n\ge 1$,
$v^\eta_{t\land \tau_n}$ is the unique fixed point of $H_{\cdot\land\tau_n}$ in $\C_{k,n}.  $

(2) We now verify \eqref{KKL}. By \eqref{PSV} and \eqref{DN1},  we find  constants $c_1,c_2>0$ such that
 \beg{align*} &\rr_{\ll}(H_{\cdot\land\tau_n}(v^1), H_{\cdot\land\tau_n} (v^2))^k= \E\bigg[\sup_{t\in [0,\tau_n]} \e^{-\ll t}|H_t(v^1)-H_t(v^2)|^k\bigg]\\
&\le c_1  \E\bigg[\sup_{t\in [0,\tau_n]} \e^{-\ll t} \bigg\{\bigg(\int_0^t \big\{ |v_s^1-v_s^2|+ \big(\E |v_s^1-v_s^2|^k\big)^{\ff 1 k}\big\}\d s\bigg)^k \\
&\qquad  \qquad \qquad \qquad\qquad +     \bigg(\int_0^t   |v_s^1-v_s^2|^2f_s(X_s )^2 \d s\bigg)^{\ff k 2} \bigg\}\bigg]\\
&\le 2c_1 T^{k-1} \rr_{\ll,n}(v^1,v^2) \sup_{t\in [0,T]} \int_0^t \e^{-\ll (t-s)}\d s\\
&\qquad + c_1 \E  \bigg[ \sup_{t\in [0,\tau_n]} \Big( \e^{-\ll t}|v_t^1-v_t^2|^k\Big) \bigg(\int_0^{t} \e^{-\ff{2\ll(t-s)}k}  f_s(X_s )^2 \d s\bigg)^{\ff k 2}\bigg]\\
&\le \rr_\ll(v^1,v^2) \bigg\{ \ff{c_2}{\ll}  +  c_1 \sup_{\OO}   \sup_{t\in [0,\tau_n]} \bigg(\int_0^{t} f_s(X_s )^{2\theta}\d s\bigg)^{\ff{k}{2\theta}} \bigg(\int_0^{t}  \e^{-\ff{2\theta^*\ll(t-s)}k}\d s\bigg)^{\ff k{2\theta^*}}\bigg\}\\
&\le  \bigg\{ \ff{c_2}{\ll}  +     c_1 n^{\ff{k}{2\theta}} \Big( \ff k{2\ll \theta^*}\Big)^{\ff k{2\theta^*}}\bigg\} \rr_\ll(v^1,v^2),\ \ v^1,v^2\in \C_{k,n}.\end{align*}
Therefore, when $\ll>0$ is large enough, $H_{\cdot\land\tau_n}$  is contractive in $\rr_{\ll}$ for large $\ll>0$, and hence has a unique fixed point on $\C_{k,n}. $

\end{proof}

 \subsection{Proofs of Theorem \ref{TA2'}(1)}

 Theorem \ref{TA2'}(1) is implied by the following result for     $\eta=\phi(X_0^\mu)$.

\beg{prp} \label{PP1} Assume $(H)$.  For any    $v\in\R^d$ and $\eta\in L^k(\OO\to\R^d,\F_0,\P)$,  $\nn_\eta X_t$ and $\nn_v X_t^{x}$ exist in $L^k(\OO\to C([0,T];\R^d),\P)$,
and for any $j\ge 1$ there exists a  constant $c>0$   such that
\beq\label{DFF-}  \E\Big[\sup_{t\in [0,T]}|\nn_\eta X_t^\mu |^j\Big|\F_0\Big]\le c \big(\{\E[|\eta|^k]\}^{\ff j k}+|\eta|^j\big),\ \  \mu\in \scr P_k,  \eta\in L^k(\OO\to\R^d,\F_0,\P),\end{equation}
\beq\label{DFF'-}  \E\Big[\sup_{t\in [0,T]}|\nn_v X_t^{\mu,x}|^j\Big]\le c |v|^j,\ \   x,v\in\R^d, \mu\in \scr P_k.\end{equation} \end{prp}

\beg{proof} The existence of $\nn_v X_t^{\mu, x}$ and \eqref{DFF'-} follow  from  Theorem \ref{T03}(1)    for  $b_t(x):= b_t(x,\mu_t)$
where the constant in  \eqref{DFF0} is uniformly in $\mu_t$ according to Remark 2.1.
 So, it suffices to prove \eqref{DFF-}.   We simply denote
 $$ X_t=X_t^\mu,\ \  v_t=v_t^\eta,\  \ \ t\in [0,T].$$
For any $r\in (0,1]$ let  $Y_t^r$ be in \eqref{YTV}. We have $Y_t:= Y_t^0= \Theta_t(X_t)$.  Let
\beq\label{YEE2} \tt v_t^\vv:=\ff{Y_t^{\vv}-Y_t }\vv,\ \ t\in [0,T], \vv\in (0,1).\end{equation}
By Theorem \ref{TA2}(2) and \eqref{AC}, for any $j\ge 1$ there exists $c(j)>0$ such that
\beq\label{TX6} \E\Big[\sup_{t\in [0,T]} |\tt v_t^\vv|^j\Big|\F_0\Big] \le  c(j)\big(\{\E[|\eta|^k]\}^{\ff j k} +  |\eta|^j)\big),\ \ \vv\in (0,1).\end{equation}
We claim that it suffices to prove
\beq\label{*N} \lim_{\vv\downarrow 0} \E\bigg[\sup_{t\in [0,T]} |\tt v_t^\vv- v_t |^k\bigg]=0.\end{equation}
Indeed, this implies that
$$\nn_\eta Y_t := \lim_{\vv\downarrow 0} \tt v_t^\vv =v_t $$ exists in $L^k(\OO\to C([0,T];\R^d); \P)$, so that \eqref{AC}, \eqref{AC''}  and $\Theta_t := id +u_t $ yield
$$\nn_\eta X_t := \lim_{\vv\downarrow 0} \ff{X_t^{\vv}-X_t }\vv  =(\nn \Theta_t (X_t ))^{-1} v_t $$ exists in the same space, and \eqref{DFF-} follows from \eqref{EST1}.

Recall that $\mu_s^\vv= \L_{X_s^\vv}, \vv\in [0,1]$.
  By \eqref{PDEM} and It\^o's formula, we  obtain
\beq\label{VEP}   \beg{split}   \d \tt v_t^\vv&=     \ff 1 \vv \Big\{\ll u(X_t^\vv)+b_t^{(1)}(X_t^\vv,\mu_t)-\ll u(X_t)  - b_t(X_t,\mu_t)\\
&\qquad \qquad +\nn_{b^{(1)}_t(X_t^{\vv},\mu_t^{\vv})- b^{(1)}_t (X_t^{\vv}, \mu_t )}\Theta_t(X_t^{\vv})\Big\}    \d t \\
 &+ \ff 1 \vv \Big\{(\nn\Theta_t ) \si_t\big\} (X_t^{\vv}) - (\nn\Theta_t ) \si_t\big\} (X_t )\Big\}  \d W_t,\ \ \tt v_0^\vv = \ff{\Theta_0 (X_0^{\vv})-\Theta_0 (X_0 )}\vv.
 \end{split}  \end{equation}
Then
 \beq\label{DN0} \beg{split}   \tt v_t^\vv=  &\,\tt v_0^\vv+\int_0^t\Big\{\nn_{(\nn\Theta_s (X_s ))^{-1} \tt v_s^{\vv}}  \{b_s^{(1)}(\cdot, \mu_s )+\ll u\} (X_s) +  \nn_{\psi_s(\tt v_s^{\vv})} \Theta_s(X_s)\Big\} \d s\\
 &+ \int_0^t \nn_{(\nn\Theta_s (X_s ))^{-1}\tt v_s^{\vv}}\{(\nn\Theta_s ) \si_s\} (X_s ) \d W_s+\aa_t^{\vv},\ \ t\in [0,T], \end{split}\end{equation}
where   $\psi_t(v)$  is in \eqref{PSV0},   and for $t\in [0,T],$
\beg{align*} &\aa_t^\vv:=\int_0^t \xi_s^\vv\d s+\int_0^t \eta_s^\vv\,\d W_s,\\
&\xi_s^{\vv}:=  \ff 1 \vv \Big\{[b^{(1)}_s(\cdot, \mu_s)+\ll u_s](X_s^{\vv})-[ b_s^{(1)}(\cdot,\mu_s) +\ll u_s](X_s) +\nn_{b^{(1)}_s(X_s^{\vv},\mu_s^{\vv})- b^{(1)}_s(X_s^{\vv}, \mu_s)}\Theta_s(X_s^{\vv})\Big\} \\
&\qquad -\Big\{\nn_{(\nn\Theta_s (X_s ))^{-1} \tt v_s^{\vv}} [ b_s^{(1)}(\cdot, \mu_s ) +\ll u_s](X_s) +  \nn_{\psi_s(\tt v_s^{\vv})} \Theta_s(X_s)\Big\},\\
&\eta_s^{\vv}:=  \ff{\{(\nn\Theta_s ) \si_s\}(X_s^{\vv})- \{(\nn\Theta_s ) \si_s\}(X_s )}\vv- \nn_{(\nn\Theta_s(X_s ))^{-1}\tt v_s^{\vv}} \{(\nn\Theta_s ) \si_s\}(X_s).\end{align*}
We claim
\beq\label{LTT} \lim_{\vv\to 0} \E\Big[\sup_{t\in [0,T]} |\aa_t^{\vv}|^n \Big|\F_0\Big]  =0,\ \ n\ge 1.\end{equation}
This can be proved by the argument leading to \eqref{DN3}, but with the conditional expectation $\E[\cdot|\F_0]$ replacing the expectation.

Firstly,  by \eqref{TX6}, $Y_t^\vv= X_t^\vv+u_t(X_t^\vv)$  and \eqref{AC}, for any $j\ge 1$ there exists $c(j)>0$ such that
 \beq\label{DXY} \sup_{\vv\in (0,1]} \E\bigg[\sup_{t\in [0,T]} \Big| \ff{X_t^{\vv}-X_t }\vv\Big|^j\Big|\F_0\Big] \le  c(j)\big(\{\E[|\eta|^k]\}^{\ff j k} +  |\eta|^j)\big).\end{equation}
 Since $\{(\nn\Theta_s ) \si_s\},   b_s^{(1)}(\cdot,\mu_s)$ and $\nn u_s $ are a.e. differentiable,  by the same reason leading to \eqref{000}, \eqref{DXY} implies that for any $s\in (0,T]$, $\P$-a.s.
  \beg{align*}  &\lim_{\vv\to 0}   \Big|\ff{\{(\nn\Theta_s ) \si_s\}(X_s^{\vv})-\{(\nn\Theta_s ) \si_s\}(X_s )}\vv-  \nn_{(\nn\Theta_s (X_s ))^{-1} \tt v_s^{\vv}} \{(\nn\Theta_s) \si_s\}(X_s )\Big|=0,\\
&\lim_{\vv\to 0} \Big|\ff{[b_s^{(1)}(\cdot, \mu_s) +u_s](X_s^{\vv})-   [b_s^{(1)} (\cdot, \mu_s)+\ll u_s](X_s)}\vv- \nn_{(\nn\Theta_s(X_s))^{-1}\tt v_s^{\vv}}   [b^{(1)}(\cdot,\mu_s)+\ll u_s](X_s)\Big| =0.\end{align*}
Next,  as in \eqref{*0*}, by the chain rule in Theorem 2.1 of \cite{BRW20} and $b^{(1)}_t\in \D_k$,   we obtain
$$  \lim_{\vv\to 0}  \Big|\ff{b_s^{(1)}(X_s^{\vv}, \mu_s^{\vv})-  b_s^{(1)}(X_s^{\vv}, \mu_s )}\vv- \psi_s(\tt v_s^{\vv}) \Big| =0,\ \ s\in (0,T].$$
Thus, for any $s\in (0,T]$, as $\vv\to 0$ we have $\P$-a.s.
\beq\label{PO1} \lim_{\vv\to 0} \big\{|\xi_s^\vv|+\|\eta_s^\vv\|\big\}=0.\end{equation}
Moreover, by $(H)$ and Lemma 2.1 in \cite{XXZZ}, we find a constant $c>0$ such that
\beg{align*}  |\xi_s^\vv|+\|\eta_s^\vv\|  \le c|\tt v_s^{\vv}| \Big(1+  \sum_{i=0}^l\big\{ \scr M f_i(s,\cdot) (X_s )+  \scr M f_i (s,\cdot)(X^{\vv}_s)\big\}\Big),\ \ s\in [0,T].\end{align*}
Finally, let $\theta>1$ be in the proof of \eqref{DN3} such that \eqref{GGR} holds for $X_t^{\vv}$ replacing $X_t^{x+\vv v}$. By \eqref{GGR} for $X_t^{\vv}$, \eqref{TX6},    and  Lemma 2.1 in \cite{XXZZ}, for any $n\ge 1$ there exist  constants $c_1(n),c_2(n)>0$ such that
\beg{align*} &\E\bigg[\bigg(\int_0^T I_s^{2\theta} \d s\bigg)^n\bigg|\F_0\bigg]\\
& \le c_1(n) \E\bigg[\Big(\sup_{s\in [0,T]} |\tt v_s^\vv|^{2\theta n}\Big) \bigg(\int_0^T \Big(1+\sum_{i=0}^l \big\{\scr M f_i^{2\theta} (s,X_s ) +\scr M f_i^{2\theta}(s,X_s^{\vv}) \big\}\Big)\d s \bigg)^{n} \bigg|\F_0\bigg]\\
& \le c_1(n) \Big(\E\Big[\sup_{s\in [0,T]} |\tt v_s^{\vv}|^{4\theta n}\Big|\F_0\Big] \Big)^{\ff 1 2}\\
&\qquad\times  \bigg(\E\bigg[\bigg(\int_0^T  \Big(1+\sum_{i=0}^l \big\{\scr M f_i^{2\theta}(s,X_s ) +\scr M f_i^{2\theta} (s,X_s^{\vv}) \big\}\Big)\d s \bigg)^{2n} \bigg|\F_0\bigg]\bigg)^{\ff 1 2} \\
&\le c_2(n) (1+|\eta|^{2\theta n})<\infty.\end{align*}
By  BDG's inequality and the dominated convergence theorem, this and \eqref{PO1} imply  \eqref{LTT}.

Now,  by  \eqref{VS} and \eqref{DN0},    the argument leading to \eqref{TTYU} gives
\beg{align*}  |v_t -\tt v_t^{\vv}|^{2k} \le &\,|v_0 - \tt v_0^\vv|^{2k}+   \int_0^t \big\{ |v_s -\tt v_s^\vv|^{2k}\gg_t + \big(\E[|v_s -\tt v_s^\vv|^k] \big)^{2}\big\}  \d t\\
&+ K  \sup_{r\in [0,t]}|\aa_r^{\vv}|^{2k} +  M_t,\ \ t\in [0,T],\end{align*}
where $K>0$ is a constant and $\gg_t$ is a positive process satisfying
$$\E \big[\e^{N \int_0^T \gg_t\d t}\big]<\infty,\ \ N>0.$$
Therefore, by the stochastic Gronwall inequality  \cite[Lemma 3.7]{XZ},  we find a constant $c>0$ such that
\beq\label{XIN} \beg{split} &\E\Big[\sup_{s\in [0,t]} |\tt v_s^\vv- v_s|^k\Big|\F_0\Big]\\
&\le c  |v_0 - \tt v_0^\vv|^{k}+ c \bigg(\E\bigg[\sup_{s\in [0,t]}|\aa_s^\vv|^{2k}\bigg|\F_0\bigg]\bigg)^{\ff 1 2}+c\bigg(\int_0^t \big(\E[\tt v_s^\vv-v_s |^k ]\big)^{2}\d s\bigg)^{\ff 1 2},\ \ t\in [0,T].\end{split}\end{equation}
By  \eqref{LTT}, \eqref{XIN} and noting that   $\lim_{\vv\to 0} |v_0-\tt v_0^\vv|=0$,  we obtain
\beq\label{MMA} \limsup_{\vv\to 0} \E\Big[\sup_{s\in [0,t]} |\tt v_s^\vv- v_s|^k\Big|\F_0\Big]\le c\limsup_{\vv\to 0}\bigg( \int_0^t \big(\E[\tt v_s^\vv-v_s |^k ]\big)^{2}\d s\bigg)^{\ff 1 2}.\end{equation}
Taking $j=k$ in \eqref{EST1} and \eqref{TX6} and \eqref{DXY} we see that
$$\bigg\{  \E\Big[\sup_{t\in [0,T]}\big\{ |\tt v_t^\vv|^k+ |v_t|^k\big\}\Big|\F_0\Big]:\ \vv\in (0,1]\bigg\}$$ is uniformly integrable with respect to $\P$, so that by Fatou's lemma, \eqref{MMA} implies
\beg{align*} h_t&:= \limsup_{\vv\to 0} \E\bigg[\sup_{s\in [0,t]}|\tt v_s^\vv- v_s|^k\bigg]=\limsup_{\vv\to 0} \E\bigg\{\E\bigg[\sup_{s\in [0,t]}|\tt v_s^\vv- v_s|^k\bigg|\F_0\bigg]\bigg\}\\
& \le \E\bigg\{\limsup_{\vv\to 0} \E\Big[\sup_{s\in [0,t]}\big\{|\tt v_s^\vv- v_s|^k\Big|\F_0\Big]\bigg\}
\le c \bigg(\int_0^t h_s^2\d s\bigg)^{\ff 1 2},\ \ t\in [0,T]\end{align*} and $h_t<\infty$,
so that $h_t=0$ for all $t\in [0,T].$ Therefore,    \eqref{*N} holds and hence   the proof is finished.   \end{proof}

\subsection{Proof  of Theorem \ref{TA2'}(2)}

 For any $\eta\in L^k(\OO\to\R^d,\F_0,\P)$, $\mu\in \scr P_k$,
and  $\vv\in [0,1]$,   let $X_t^\vv$ solve \eqref{E10} for $X_0^\vv= X_0^\mu+\vv\eta$. Consider
$$ \GG_\eta(f(X_t^\mu)):=\lim_{\vv\downarrow 0} \ff{\E[f(X_t^\vv)-f(X_t^\mu)]}\vv,\ \ t\in (0,T], f\in \B_b(\R^d).$$
Theorem \ref{TA2'}(2)  is implied by the following result for  $\eta=\phi(X_0)$.

\beg{prp}\label{PP2}  Assume $(H)$.    For any $\eta\in L^k(\OO\to\R^d,\F_0,\P)$ and  $\mu\in \scr P_k$,   $D_\eta^I  P_tf(\mu) $ exists and satisfies the following formula for any    $\bb \in C^1([0,t])$ with $\bb_0=0$ and $\bb_t=1:$
\beq\label{BSM2} \beg{split} & \GG_\eta(f(X_t^\mu))
= \int_{\R^d\times \R^d} \E\bigg[f(X_t^{\mu,x}) \int_0^t \bb_s'\big\<\zeta_s(X_s^{\mu,x} ) \nn_{v}X_s^{\mu,x},  \d W_s\big\> \bigg]\L_{(X_0^\mu, \eta)}(\d x,\d v)\\
&+  \E\bigg[f(X_t^{\mu}) \int_0^t  \Big\<\zeta_s(X_s^\mu)  \E\big[\<D^Lb^{(1)}_s(z, P_s^*\mu)(X_s^\mu), \nn_\eta X_s^\mu\>\big]\big|_{z=X_s^\mu},  \d W_s\Big\> \bigg].
\end{split} \end{equation}
Consequently,     there exists a constant $c >0$ such that
\beq\label{GRDI2}\beg{split} & \big|\GG_\eta(f(X_t^\mu))\big|  \le \ff {c}{\ss t} \big(P_t |f|^{k^*}(\mu)\big)^{\ff 1 {k^*}}(\E[|\eta|^k])^{\ff 1 k},\\
&t\in (0,T], f\in \B_b(\R^d),\mu\in \scr P_k, \eta\in  L^k(\OO\to\R^d,\F_0,\P).\end{split} \end{equation}
\end{prp}
\beg{proof}  Let $X_t^{\mu,x}$ solve \eqref{E102}. Since $X_t^\mu$ solve \eqref{E102} with inital value $X_0^\mu$,    the strong uniqueness  implies
\beq\label{XXU} X_t^\mu = X_t^{\mu, X_0^\mu},\ \  t\in [0,T].\end{equation}
Let
$(P_{s,t}^\mu)_{0\le s\le t\le T}$ be the semigroup associated with    \eqref{E102}, i.e. for  $(X_{s,t}^{\mu,x})_{t\in [s,T]}$ solving \eqref{E102} from time $s$ with $X_{s,s}^{\mu,x}=x$,
\beq\label{TSM} P_{s,t}^\mu f(x):= \E[f(X_{s,t}^{\mu,x})],\ \ t\in [s,T], x\in\R^d.\end{equation}
Simply denote  $P_t^\mu= P_{0,t}^\mu$. Then
  \eqref{XXU} implies
\beq\label{EN0} P_tf(\mu)= \E[f(X_t^\mu)]= \int_{\R^d} P_t^\mu f(x)\mu(\d x),\ \ t\in [0,T], f\in \B_b(\R^d).\end{equation}
By Theorem \ref{T03}, $(H)$ implies that     for any $t\in (0,T]$ and $\bb\in C^1([0,t])$ with $\bb_0=0$ and $\bb_t=1$,
 \beq\label{EN1} \nn_v P_t^\mu f(x)= \E\bigg[f(X_t^{\mu,x}) \int_0^t \bb_s'\big\<\zeta_s(X_s^{\mu,x})\nn_v X_s^{\mu,x}, \d W_s\>\bigg],\ \ v\in \R^d, f\in\B_b(\R^d).\end{equation}
 Next, denote $\mu_t=P_t^*\mu=\L_{X_t^\mu}$ and  let ${\bar X}_s^{\vv}$ solve \eqref{E102} for ${\bar X}_0^\vv= X_0^\vv$, i.e.
 \beq\label{*EN} \d \bar X_s^\vv = b_s(\bar X_s^\vv,\mu_s)\d s +\si_s(\bar X_s^\vv)\d W_s,\ \  s\in [0,t], \bar X_0^\vv= X_0^\vv.\end{equation}
We have
\beg{align*} &\E[ f(\bar X_t^\vv)]= \int_{\R^d} (P_t^\mu)(x) ]\L_{X_0^\mu+\vv \eta}(\d x)\\
&= \int_{\R^d\times\R^d} P_t^\mu f(x+\vv v)\L_{(X_0^\mu,\eta)}(\d x, \d v),\ \ f\in \B_b(\R^d),\end{align*}
where $\L_{(X_0^\mu,\eta)}$ is the law of $(X_0^\mu,\eta)$.
 Combining this with \eqref{EN0} and \eqref{EN1}, and applying the dominated convergence theorem, we obtain
 \beq\label{EN2}\beg{split} & \lim_{\vv\to 0}  \ff{\E[f(\bar X_t^\vv)] - P_t f(\mu)}\vv =  \int_{\R^d\times\R^d} \nn_{v} P_t^\mu f(x)\L_{(X_0^\mu,\eta)}(\d x, \d v)\\
 &= \int_{\R^d\times\R^d} \E\bigg[f(X_t^{\mu,x}) \int_0^t \bb_s' \<\zeta_s(X_s^{\mu,x}) \nn_{v} X_s^{\mu,x}, \d W_s\>\bigg]\L_{(X_0^\mu,\eta)}(\d x, \d v).\end{split}\end{equation}

 On the other hand, denote $\mu_t^\vv= \L_{X_t^\vv}$ and   let
 $$R_t^\vv:= \e^{\int_0^t\<\zeta_s(X_s^\vv) \{b_s^{(1)}(X_s^\vv,\mu_s)- b_s^{(1)}(X_s^\vv,\mu_s^\vv)\}, \d W_s\>-\ff 1 2\int_0^t |\zeta_s(X_s^\vv) \{b_s^{(1)}(X_s^\vv,\mu_s)- b_s^{(1)}(X_s^\vv,\mu_s^\vv)\}|^2\d s}.$$
By $(H)$, $\zeta_s= \si_s^*(\si_s\si_s^*)^{-1}$ and Girsanov's theorem, $\Q_t^\vv:= R_t^\vv\P$ is a probability measure under which
$$\tt W_r^\vv:=W_r- \int_0^r\zeta_s(X_s^\vv) \{b_s^{(1)}(X_s^\vv,\mu_s)- b_s^{(1)}(X_s^\vv,\mu_s^\vv)\}\d s,\ \ r\in [0,t]$$ is a Brownian motion, and
\beq\label{LLN} \sup_{r\in [0,T], \vv\in (0,1]} \E\Big[\ff{|R_r^\vv-1|^j}{\vv^j}\Big]<\infty,\ \ j\ge 1. \end{equation}
Reformulate the SDE  for $X_s^\vv$ as
$$ \d X_s^\vv = b_s(X_s^\vv, \mu_s) +\si_s(X_s^\vv)\d \tt W_s^\vv,\ \ X_0^\vv= \bar X_0^\vv.$$
By the well-posedness we obtain
$\L_{X_t^\vv|\Q_t^\vv} = \L_{\bar X_t^\vv|\P}$, so that
$$ \E[f(\bar X_t^\vv)]= \E[R_t^\vv f(X_t^\vv)],\ \ f\in \B_b(\R^d).$$
Thus,
 \beg{align*} & \ff{\E[ f(X_t^\vv)]- \E[f(\bar X_t^\vv)]}\vv =  \ff{\E[f(X_t^\vv)(1-R_t^\vv)]}\vv= I_1(\vv)+ I_2(\vv),\\
& I_1(\vv):= \E\bigg[f(X_t^\mu) \ff{1-R_t^\vv}\vv\bigg],\ \ I_2(\vv):= \E\bigg[\{f(X_t^\vv)-f(X_t^\mu)\} \ff{1-R_t^\vv}\vv\bigg].\end{align*}
By  \eqref{BKK},   \eqref{*0*} and the dominated convergence theorem,  we obtain
$$\lim_{\vv\to 0} I_1(\vv)
 = \E\bigg[ f(X_t^\mu) \int_0^t \big\<\zeta_s(X_s^\mu)  \E[\<D^L b_s^{(1)}(z,\mu_s), \nn_\eta X_s^\mu\>]|_{z=X_s^\mu},  \d W_s\big\>\bigg].$$
So, to prove  \eqref{BSM2} it suffices to verify
\beq\label{LN10} \lim_{\vv\to 0} I_2(\vv)=0.\end{equation}
 By \eqref{BKK}, we have
\beq\label{LN11} \lim_{r\uparrow t} \sup_{\vv\in (0,1]} \E\Big[\ff{|R_t^\vv-R_r^\vv|}\vv\Big]=0.\end{equation}
Since $(H)$ holds for $[r,T]$ replacing $[0,T]$, \eqref{GRD00} holds for $(r,T)$ replacing $(0,T)$. Similarly,
\eqref{GRD} holds for $P_{r,t}^{\mu^\vv}$ and $P_{r,t}^\mu$ defined in \eqref{TSM}  replacing $P_{t-r}$. Therefore, by the Markov property,
\beq\label{NL0} \beg{split}& |\E[f(X_t^\vv)-f(X_t^\mu)|\F_{r}]|= |(P_{r,t}^{\mu^\vv} f)(X_{r}^\vv)- (P_{r, t}^\mu f)(X_{r}^\mu)| \\
&\le  |(P_{r,t}^{\mu^\vv} f)(X_{r}^\vv)- (P_{r, t}^{\mu^\vv} f)(X_{r}^\mu)|+  |(P_{r,t}^{\mu^\vv} f)(X_{r}^\mu)- (P_{r, t}^\mu f)(X_{r}^\mu)|\\
&\le  c \|f\|_\infty\bigg(\ff{|X_{r}^\vv- X_{r}^\mu|}{\ss{t-s}} \land 1\bigg) +  |(P_{r,t}^{\mu^\vv} f)(X_{r}^\mu)- (P_{r, t}^\mu f)(X_{r}^\mu)|.\end{split}\end{equation}
On the other hand, let $(\tt X_{r,s}^\vv)_{s\in [r,t]}$ solve the SDE
$$\d \tt X_{r,s}^\vv = b_s(\tt X_{r,s}^\vv, \mu_s^\vv)\d s +\si_s(\tt X_{r,s}^\vv)\d W_s,\ \ \tt X_{r,r}^\vv= X_{r}^\mu, s\in [r,t].$$
We have
$$P_{r,t}^{\mu^\vv} f(X_{r}^\mu)= \E\big[f(\tt X_{r,t}^\vv)\big|\F_{r}\big],\ \ P_{r,t}^\mu f(X_{r}^\mu)= \E\big[f(X_t^\mu)|\F_{r}\big].$$
Noting that \eqref{BKK} and \eqref{LN2} imply
\beq\label{NL2} |b(x,\mu_t^\vv)- b_t(x,\mu_t)|\le c_1 \W_k(\mu_0^\vv,\mu_0)\le c_1 \vv (\E[|\eta|^k])^{\ff 1 k}  \end{equation}
for some constant $c_1>0$, by Girsanov's theorem,
$$R_{r,t}^\vv:=\e^{\int_{r}^t \<\zeta(X_s^\mu) \{b_s(X_{s}^\mu,\mu_s^\vv)- b_s(X_s^\mu,\mu_s)\},\d W_s\>-\ff 1 2\int_{r}^t |\zeta(X_s) \{b_s(X_{s}^\mu,\mu_s^\vv)- b_s(X_s^\mu,\mu_s)\}|^2\d s}$$
is a probability density such that under $\Q_{r,t}:= R_{r,t}^\vv \P$,
$$\tt W_s:= W_s- \int_{r}^s \zeta(X_\theta) \{b_s(X_{\theta}^\mu,\mu_\theta^\vv)- b_\theta(X_\theta^\mu,\mu_\theta)\}\d \theta,\ \ s\in [r,t]$$
is a Brownian motion. Reformulating the SDE for $(X_{s}^\mu)_{s\in [r,t]}$ as
$$\d X_s^\mu= b_s(X_s^\mu, \mu_s^\vv)\d s +\si_s(X_s^\mu) \d \tt W_s,\ \ X_{r}^\mu= \tt X_{r,r}^\vv,\ \ s\in [r,t], $$
by the uniqueness we obtain
$$P_{r,t}^{\mu^\vv} f(X_{r}^\mu)= \E\big[R_{r,t}^\vv f(X_t^\mu)\big|\F_{r}\big],$$
so that by   Pinsker's inequality  and \eqref{NL2}, we find a  constant  $c_2>0$ such that
\beq\label{NL*} \beg{split} &|(P_{r,t}^{\mu^\vv} f)(X_{r}^\mu)- (P_{r, t}^\mu f)(X_{r}^\mu)|^2\le \|f\|_\infty \big|\E[|1-R_{r,t}^\vv|\big|\F_{r}]\big|^2\\
&\le 2\|f\|_\infty\E_{\Q_{r,t}}\big[ \log R_{r,t}^\vv\big|\F_0\big] \\
&= \|f\|_\infty\int_{r}^t \E_{\Q_{r,t}}\big[|\zeta(X_s^\mu) \{b_s(X_{s}^\mu,\mu_s^\vv)- b_s(X_s^\mu,\mu_s)\}|^2\big| \F_{r}\big] \d s\\
&\le   c_2 \|f\|_\infty(t-r)\vv^2\|\eta\|_{L^k(\P)}^2.\end{split}\end{equation}
Combining this with \eqref{EE0*}, \eqref{LLN} and  \eqref{NL0}, and noting that $(s\land 1)^2\le s$ for $s\ge 0$, we find constants $c_3,c_4>0$ such that
\beg{align*}&  \bigg|\E\Big[\{f(X_t^\vv)-f(X_t)\}\ff{1-R_r^\vv}\vv\Big] \bigg|\le \bigg(\E\Big|\E\big[f(X_t^\vv)-f(X_t)\big|\F_r\big]\Big|^2\bigg)^{\ff 1 2} \bigg(\E\Big[\ff{|1-R_r^\vv|^2}{\vv^2} \Big]\bigg)^{\ff 1 2}\\
&\le  c_4\|f\|_\infty\bigg(  \ff{\E[|X_{r}^\vv-X_{r}^\mu|]}{\ss{t-r}}  \bigg)^{\ff 1 2}   + c_4\|f\|_\infty\vv\\
 & \le c_5\ss T\|f\|_\infty\Big(\ff{\vv}{t-r}\Big)^{\ff 1 2},\ \ \vv\in (0,1], t\in [0,T].\end{align*}
Combining this with \eqref{LN11} we obtain
$$\lim_{\vv\downarrow  0} I_2(\vv)\le \lim_{r\uparrow  t} \lim_{\vv\downarrow 0} \bigg\{\bigg|\E\Big[\{f(X_t^\vv)-f(X_t)\}\ff{1-R_r^\vv}\vv\Big]\bigg|+2 \|f\|_\infty \E\Big[\ff{|R_t^\vv-R_r^\vv|}\vv\Big]\bigg\} =0. $$
Therefore, \eqref{LN11} holds.

It remains to   prove  \eqref{GRDI2}. By Jensen's inequality, it suffices to prove for $p\in (1, 2]$. By \eqref{BSM2}, we have
\beq\label{AA0} |\GG_\eta(f(X_t^\mu))|\le \E(|J_1(X_0^\mu,\eta)|) +| J_2|,\end{equation} where
\beg{align*} &J_1(x,v):=   \E\bigg[f(X_t^{\mu,x}) \int_0^t \bb_s'\big\<\zeta_s(X_s^{\mu,x} ) \nn_{v}X_s^{\mu,x},  \d W_s\big\> \bigg], \\
&J_2:=   \E\bigg[f(X_t^{\mu}) \int_0^t  \Big\<\zeta_s(X_s^\mu)  \E\big[\<D^Lb^{(1)}_s(z, P_s^*\mu)(X_s^\mu), \nn_\eta X_s^\mu\>\big]|_{z=X_s^\mu},  \d W_s\Big\> \bigg].\end{align*}
Taking $\bb_s=\ff s t$, by $\|\zeta\|_\infty<\infty$,  \eqref{DFF'} and H\"older's inequality, we find   constants $c_1, c_2>0$ such that
\beg{align*}   |J_1(x,v)|  &\le \ff{c_1}t  \big(P_t^\mu |f|^{p}(x)\big)^{\ff 1 {p}} \bigg\{\E\bigg[\bigg( \int_0^t  |\nn_v X_s^{\mu,x}|^2\d s\bigg)^{\ff {p^*}2}\bigg]\bigg\}^{\ff 1{p^*}}\\
&\le \ff {c_2|v|}{\ss t}  (P_t^\mu |f|^{p}(x))^{\ff 1{p}},\ \ t\in (0,T]. \end{align*}
Combining this with \eqref{EN0} and $P_t^\mu|f|^{p}(X_0^\mu)=\E[|f(X_t^\mu)|^p|\F_0]$, we derive
\beq\label{AA1}\beg{split} & \E[|J_1(X_0^\mu,\eta)|]\le   \ff {c_2}{\ss t}  \E\Big[|\eta|\big(P_t^\mu|f|^{p}(X_0^\mu)\big)^{\ff 1 {p}}\Big]\\
& \le \ff {c_2 \|\eta\|_{L^k(\P)}}{\ss t}  \big\|\big(\E[|f(X_t^\mu)|^p|\F_0]\big)^{\ff 1 p}\big\|_{L^{k^*}(\P)},\ \ t\in (0,T].\end{split}\end{equation}
On the other hand, by $(H)$, H\"older's inequality  and $\eqref{DFF-}$ for $j=k$,  we find   constants  $c_3, c_4>0$ such that
$$ I_s(z):= \Big|\zeta_s(X_s^\mu)  \E\big[\<D^Lb^{(1)}_s(z, P_s^*\mu)(X_s^\mu), \nn_\eta X_s^\mu\>\big]\Big| \le c_3 \| \nn_\eta X_s^\mu\|_{L^k(\P)} \le c_4 \|\eta\|_{L^k(\P)},$$
so that
\beg{align*}   |J_2| & \le \E\bigg[\big(\E[|f(X_t^\mu)|^p\F_0]\big)^{\ff 1 p} \bigg(\E \bigg[ \int_0^t I_s(X_s^\mu)^2\d s\bigg]^{\ff{p^*} 2} \bigg)^{\ff 1 {p^*}} \bigg]\\
& \le c_4 \ss {  t} \|\eta\|_{L^k(\P)} \E\Big[\big(\E[|f(X_t^\mu)|^p\F_0]\big)^{\ff 1 p}\Big].\end{align*}
This and \eqref{AA1} imply  \eqref{GRDI2}.
\end{proof}

\subsection{Proof  of Theorem \ref{TA2''}}

Simply denote $X_t=X_t^\mu,$ and for any $\vv\in [0,1]$ let
$X_t^\vv$ solve \eqref{E10} with $X_0^\vv= X_0+\vv \phi(X_0),$  $\mu^\vv:=\L_{X_0+\vv\phi(X_0)}$ and $\mu_t^\vv:= P_t^*\mu^\vv=\L_{X_t^\vv}.$ We have
$$P_tf(\mu\circ (id+\vv\phi^{-1}))=\E[f(X_t^\vv)].$$ It suffices to prove
\beq\label{LMT} \lim_{\vv\downarrow 0} \sup_{\|\phi\|_{L^k(\mu)}\le 1} \Bigg| \ff{ \E [f(X_t^\vv) - f(X_t)]}\vv- D_\phi^I P_tf(\mu)\bigg| =0.\end{equation}
By   applying  \eqref{BSM2} with $\bb_s=\ff s t$ for $(\mu^r, \phi(X_0))$ replacing $(\mu,\eta)$,  we obtain
\beg{align*} &\ff{\d}{\d r} \E [f(X_t^r)]:= \lim_{\vv\downarrow 0} \ff{\E[f(X_t^{r+\vv})- f(X_t^r)]}\vv=\GG_{\phi(X_0)}(f(X_t^{\mu^r}))\\
& =\ff 1 t  \int_{\R^d} \E\bigg[f(X_t^{\mu^r,x+r\phi(x)}) \int_0^t  \big\<\zeta_s(X_s^{\mu^r,x+r\phi(x)} ) \nn_{\phi (x)}X_s^{\mu^r,x+r\phi(x)},  \d W_s\big\> \bigg]\mu(\d x)\\
&\quad +  \E\bigg[f(X_t^{r}) \int_0^t  \Big\<\zeta_s(X_s^r)  \E\big[\<D^Lb^{(1)}_s(z, \mu_s^r)(X_s^r), \nn_{\phi (X_0)} X_s^{\mu^r}\>\big]|_{z=X_s^r},  \d W_s\Big\> \bigg]. \end{align*}
Combining this  with \eqref{BSMI} for $\bb_s=\ff s t$, we derive
\beg{align*} &\sup_{\|\phi\|_{L^k(\mu)}\le 1}\bigg| \ff{ \E [f(X_t^\vv) - f(X_t)]}\vv- D_\phi^I P_tf(\mu)\bigg| \\
&= \sup_{\|\phi\|_{L^k(\mu)}\le 1}\bigg| \ff 1 \vv\int_0^\vv \Big\{\ff{\d}{\d r} \E[f(X_t^r)] - D_\phi^I P_t f(\mu)\Big\}\d r\bigg| \le \ff c{t\vv}\int_0^\vv \sum_{i=1}^4 \aa_i(r)\d r\end{align*}
for some constant $c>0$, where letting
\beg{align*} &F_\phi(r,x):=  \int_0^t \big\<\zeta_s(X_s^{\mu^r,x+r\phi(x)} ) \nn_{\phi (x)}X_s^{\mu^r,x+r\phi(x)},  \d W_s\big\>,\ \ r\in [0,1], x\in\R^d, \\
&G_\phi(r):=   \int_0^t  \Big\<\zeta_s(X_s^r)  \E\big[\<D^Lb^{(1)}_s(z, \mu_s^r)(X_s^r), \nn_{\phi(X_0)} X_s^{\mu^r}\>\big]|_{z=X_s^r},\d W_s\Big\>,\ \ r\in [0,1],\end{align*}  we set
\beg{align*}& \aa_1(r):= \sup_{\|\phi\|_{L^k(\mu)}\le 1}\bigg| \int_{\R^d} \E\Big[ f(X_t^{\mu^r, x+r\phi(x)})  \big\{F_\phi(r,x)- F_\phi(0,x)\big\}\Big]\mu(\d x)\bigg|,\\
&\aa_2(r):=  \sup_{\|\phi\|_{L^k(\mu)}\le 1}\bigg|\int_{\R^d} \E\Big[ \big\{f(X_t^{\mu^r, x+r\phi(x)}) -f(X_t^{\mu,x})\big\} F_\phi(0,x)\Big]\mu(\d x)\bigg|,\\
& \aa_3(r):=  \sup_{\|\phi\|_{L^k(\mu)}\le 1}\Big|\E\Big[f(X_t^{r})\big\{G_\phi(r)-G_\phi(0)\big\}  \Big]\Big|,\\
& \aa_4(r):= \sup_{\|\phi\|_{L^k(\mu)}\le 1}\Big|\E\Big[\big\{f(X_t^{r})-f(X_t)\big\} G_\phi(0) \Big]\Big|.\end{align*}
Since $\|f\|_\infty<\infty$, by $(H)$, \eqref{DFF} and \eqref{DFF'}, we conclude that $\{\aa_i\}_{1\le i\le 4}$ are bounded on $[0,1]$. So, \eqref{LMT} follows if
$$\lim_{r\downarrow 0} \aa_i(r)=0,\ \ 1\le i\le 4.$$
 To prove these limits, we need   the following two lemmas.

\beg{lem}\label{LM1} Assume $(H)$. For any $j\ge 1$ there exists a constant $c>0$ such that for any $\mu\in \scr P_k$ and $\phi\in T_{\mu,k}$ with $\|\phi\|_{L^k(\mu)}\le 1$,
$$  \E\Big[\sup_{t\in [0,T]}  |X_t^{\mu^r, x+r\phi(x)}- X_t^{\mu,x}|^j\Big]\le c r^j(1+|\phi(x)|^j),\ \ r\in [0,1].$$
\end{lem}

 \beg{proof} By \eqref{DFF'},  we have
 $$  \E[ |X_t^{\mu^r, x+r\phi(x)}- X_t^{\mu^r,x}|^j]\le c r^j|\phi(x)|^j,\ \ r\in [0,1], x\in\R^d.$$
Combining this with  $\W_k(\mu^r,\mu)\le r \|\phi\|_{L^k(\mu)}\le r$, we need only
 to prove
 \beq\label{LM11}\sup_{x\in\R^d}  \E[ |X_t^{\mu, x}- X_t^{\nu,x}|^j]\le c \W_k(\mu,\nu)^j,\ \  \mu,\nu\in\scr P_k\end{equation}
 for some constant $c>0$, where $X_t^{\nu,x}$ solves \eqref{E102} for $\nu_t:= P_t^*\nu$ replacing $\mu_t:=P_t^*\mu$.
 Let $u$ solve \eqref{PDEM} such that \eqref{AC} holds. Let $\Theta_t=id +u_t$ and
 $$Y_t^{\mu,x}:= \Theta_t(X_t^{\mu,x}),\ \ Y_t^{\nu,x}:=\Theta_t(X_t^{\nu,x}),\ \ t\in [0,T].$$
 By It\^o's formula we obtain
 \beg{align*}  \d (Y_t^{\mu,x}- Y_t^{\nu,x})= &\,\Big\<\big\{(\nn\Theta_t)\si_t\big\}(X_t^{\mu,x}) - (\nn\Theta_t)\si_t\big\}(X_t^{\nu,x}), \d W_t\Big\>\\
 &+ \Big\{b_t^{(1)}(X_t^{\mu,x},\mu_t) +\ll u_t (X_t^{\mu,x})- b_t^{(1)}(X_t^{\nu,x},\nu_t) -\ll u_t (X_t^{\nu,x})\\
 &\quad + \nn_{b_t^{(1)}(X_t^{\mu,x},\mu_t)- b_t^{(1)}(X_t^{\nu,x}, \nu_t)} u_t(X_t^{\nu,x})\Big\}\d t.\end{align*}
 By $(H)$, \eqref{AC},  Lemma 2.1 in \cite{XXZZ}   and It\^o's formula, for any $j\ge 1$ we find a constant $c>0$ such that
\beg{align*} |Y_t^{\mu,x}-Y_t^{\nu,x}|^{2j}\le &\,c \int_0^t |Y_s^{\mu,x}-Y_s^{\nu,x}|^{2j} \sum_{i=0}^l \big\{1+ \scr M f_i^2(x, X_s^{\mu,x})+  \scr M f_i^2(x, X_s^{\nu,x})\big\}\d s\\
 &+ c\int_0^t  \W_k(\mu_s,\nu_s)^{2j}\d s+ M_t,\ \ t\in [0,T]\end{align*}
 holds for some local martingale $M_t$ with $M_0=0$.  Since $\W_k(\mu_s,\nu_s)\le c\W_k(\mu,\nu)$ due to \eqref{LN2},
     \eqref{LM11} follows from the stochastic Gronwall inequality, Lemma 2.1 in \cite{XXZZ},  and Khasminskii's estimate \eqref{KHA} for $X_s^{\mu,x}$ and
 $X_s^{\nu,x}$ replacing $X_s$.
 \end{proof}

 \beg{lem}\label{LM2} Assume $(H)$ and $\eqref{00*}$. For any $j\ge 1$ there exist a constant $c>0$ and   a positive function $\vv (\cdot)$ on $[0,1]$ with $\vv(r)\downarrow 0$ as $r\downarrow 0$, such that for any $\phi\in T_{\mu,k}$
 with $\|\phi\|_{L^k(\mu)}\le 1$ and $r\in [0,1]$,
 \beq\label{LA1} \sup_{|v|\le 1} \E\bigg[\sup_{t\in [0,T]} \big|\nn_v X_t^{\mu^r, x+r\phi(x)}- \nn_v X_t^{\mu,x}\big|^{j}\le\min\big\{c,  \vv (r) (1+ |\phi(x)|^j)\big\},\ \   x\in\R^d,\end{equation}
 \beq\label{LA2}  \E\bigg[\sup_{t\in [0,T]} \big|\nn_{\phi(X_0)} X_t^{\mu^r}- \nn_{\phi(X_0)}X_t^{\mu}\big|^{j}\Big|\F_0\Big]\le |\phi(X_0)|^j\min\big\{c, \vv(r) (1+ |\phi(X_0)|^j)\big\}.\end{equation}
\end{lem}

 \beg{proof} We only prove \eqref{LA1} since \eqref{LA2} can be proved in the same way by using \eqref{DFF-} and \eqref{EE0} replacing \eqref{DFF'} and Lemma \ref{LM1} respectively.
 We simply denote
\beq\label{N01} X_t^x:= X_t^{\mu,x},\ \ X_t^{r,x}:= X_t^{\mu^r, x+r\phi(x)},\ \ \tt v_t:= \nn_v X_t^{\mu,x},\ \ \tt v_t^r:= \nn_v X_t^{\mu^r, x+r\phi(x)}.\end{equation}
Let $u$ solve \eqref{PDEM} such that \eqref{AC} holds. We may also assume that $u$ satisfies \eqref{AC''} as explained before. Let $\Theta_t= id + u_t$ and denote
\beq\label{N02} Y_t^x:= \Theta_t(X_t^x),\ \ Y_t^{r,x}:= \Theta_t(X_t^{r, x}),\ \ v_t:= (\nn \Theta_t (X_t^x))^{-1}\tt v_t,\ \ v_t^r:=  (\nn \Theta_t (X_t^{r,x}))^{-1}\tt v_t^r.\end{equation}
 By \eqref{DFF'} and \eqref{AC}, to prove \eqref{LA1}  it suffices to find $\vv(r)\downarrow 0$ as $r\downarrow 0$ such that
  \beq\label{LA1'} \sup_{|v|\le 1} \E\bigg[\sup_{t\in [0,T]} \big|v_t^r- v_t\big|^{j}\bigg]\le \vv (r) (1+ |\phi(x)|^j),\ \ r\in [0,1], x\in\R^d.\end{equation}
 By Jensen's inequality, we only need to prove for $j\ge 4$.

 To calculate $v_t$ and $v_t^r$, for any $\vv\in [0,1]$ we let
 $$Y_t^{r, x}(\vv):= \Theta_t(X_t^{\mu^r, x+r\phi(x)+\vv v}),\ \ Y_t^{x}(\vv):=\Theta_t(X_t^{\mu,x+\vv v}).$$
 Then the argument leading to \eqref{DRR} implies that
 \beq\label{N03} v_t= \lim_{\vv\downarrow  0} \ff{ Y_t^{x}(\vv)- Y_t^{x}}\vv,\ \ v_t^r=\lim_{\vv\downarrow 0} \ff{Y_t^{r,x}(\vv)- Y_t^{r,x}}\vv.\end{equation}
 By \eqref{PDEM} and It\^o's formula, we obtain
\beg{align*}  \d Y_t^{x}(\vv)=&\, \Big\{b_t^{(1)}(X_t^{\mu, x+ \vv v},\mu_t) + \ll u_t(X_t^{\mu, x+\vv v})\Big\}\d t+ \big\{(\nn\Theta_t)\si_t\big\}(X_t^{\mu, x+\vv v}) \d W_t,\\
&   \ Y_0^{x}(\vv)= x+\vv x,\end{align*}
 \beg{align*} \d Y_t^{r,x}(\vv)=&\, \Big\{b_t^{(1)}(X_t^{\mu^r, x+r\phi(x)+\vv v},\mu_t^r) + \ll u_t(X_t^{\mu^r, x+r\phi(x)+\vv v})\\
 & + \nn_{b^{(1)}_t(X_t^{\mu^r, x+r\phi(x)+\vv v},\mu_t^r) - b^{(1)}_t(X_t^{\mu^r, x+r\phi(x)+\vv v},\mu_t)}u_t(X_t^{\mu^r, x+r\phi(x)+\vv v})\Big\}\d t\\
 &+ \big\{(\nn\Theta_t)\si_t\big\}(X_t^{\mu^r,x+r\phi(x)+\vv v}) \d W_t,\ \ Y_0^{r,x}(\vv)= x+r\phi(x)+\vv x.\end{align*}
 Combining this with  \eqref{N01} and \eqref{N03}, we conclude that $v_t$ and $v_t^r$ solves the SDEs
 \beg{align*}&\d v_t= \Big\{\nn_{\tt v_t} b_t^{(1)}(X_t^x,\mu_t) +\ll \nn_{\tt v_t} u_t(X_t^x)\Big\}\d t+ \nn_{\tt v_t} \big\{(\nn\Theta_t)\si_t\big\}(X_t^x) \d W_t,\\
 & \ v_0= (\nn\Theta_0(x))^{-1} v,\end{align*}
 \beg{align*} \d v_t^r=&\, \Big\{\nn_{\tt v_t^r} b_t^{(1)}(X_t^{r,x}),\mu_t^r) + \ll \nn_{\tt v_t^r} u_t(X_t^{r,x})
   + \nn_{\tt v_t^r-\tt v_t} u_t(X_t^{r,x})  \Big\}\d t\\
 &+ \nn_{\tt v_t^r} \big\{(\nn\Theta_t)\si_t\big\}(X_t^{r,x}) \d W_t,\ \ v_0^{r}= (\nn\Theta_0(x+r\phi(x))^{-1} v.\end{align*}
Therefore, by \eqref{N02},
$$z_t^r:= v_t^r-v_t,\ \ t\in [0,T]$$ solves the SDE
\beq\label{NPP}\beg{split}   \d z_t^r =  &\,\Big\{\nn_{(\nn\Theta_t(X_t^x))^{-1} z_t^r} \big[b_t^{(1)}(\cdot,\mu_t)+\ll u_t\big](X_t^x) +\nn_{(\nn\Theta_t(X_t^x))^{-1} z_t^r} u_t(X_t^{r,x})\Big\}\d t \\
&+ \nn_{(\nn\Theta_t(X_t^x))^{-1} z_t^r} \big\{(\nn\Theta_t)\si_t\big\}(X_t^x) \d W_t - \eta_t^r\d t -\xi_t^r \d W_t,\\
& \ \ \ \qquad z_0^r=\big\{(\nn\Theta_0(x+r\phi(x))^{-1}- (\nn\Theta_0)(x))^{-1}\big\}v,\end{split}\end{equation}
where
\beg{align*}& \eta_t^r:=   \nn_{(\Theta_t(X_t^{r,x}))^{-1} v_t^r } b_t^{(1)}(X_t^x, \mu_t) - \nn_{(\Theta_t(X_t^x))^{-1} v_t^r} b_t^{(1)}(X_t^{r,x}, \mu_t^r) -\ll \nn_{(\Theta_t(X_t^{r,x}))^{-1} v_t^r } u_t(X_t^{r,x})\\
&\qquad + \ll \nn_{(\Theta_t(X_t^x))^{-1} v_t^r } u_t(X_t^x)+ \nn_{\{(\nn\Theta_t(X_t^x))^{-1}- (\nn\Theta_t(X_t^{r,x}))^{-1}\} v_t^r} u_t(X_t^{r,x}),\\
&\xi_t^r:= \nn_{(\Theta_t(X_t^x))^{-1} v_t^r } \big\{(\nn\Theta_t)\si_t\big\}(X_t^x) -  \nn_{(\Theta_t(X_t^{r,x}))^{-1} v_t^r} \big\{(\nn\Theta_t)\si_t\big\}(X_t^{r,x}),\ \ t\in [0,T].\end{align*}
By \eqref{AC},  $(H)$ and Lemma 2.1 in \cite{XXZZ}, we find a constant $c_1>0$ such that
\beg{align*} |\eta_t^r|+\|\xi_t^r\|\le &\, c_1 |v_t^r|\Big\{\|\nn b_t^{(1)}(X_t^x, \mu_t)-\nn b_t^{(1)}(X_t^{r,x},\mu_t^r)\|\\
&\qquad + |X_t^{r,x}-X_t^x| \sum_{i=0}^l \big(1+\scr M f_i(t, X_t^x) +\scr M f_i(t,X_t^{r,x})\big)\Big\}.\end{align*} By the boundedness of $\nn b^{(1)}$ and \eqref{00*}, we have
\beq\label{N0N} \|\nn b_t^{(1)}(X_t^x, \mu_t)-\nn b_t^{(1)}(X_t^{r,x},\mu_t^{r})\|\le n\big\{|X_t^x-X_t^{r,x}|+\W_k(\mu_t,\mu_t^r)\big\}^{\ff 1 {2j}} +s_n,\ \ n\ge 1,\end{equation}
where for $\varphi(r):=  \sup_{|x-x'|+ \W_k(\mu,\nu)\le r} \|\nn b_t^{(1)}(x, \mu)-\nn b_t^{(1)}(x',\nu)\|,$
$$s_n:= \sup_{r\ge 0} \big\{\varphi(r)- n r^{\ff 1 {2j}}  \big\}\downarrow 0\ \text{as}\ n\uparrow \infty.$$
Using the notation \eqref{N01}, by combining this with   Lemma \ref{LM1}, \eqref{LN2} and  \eqref{GGR} for the processes $X_t^x$ and $X_t^{r,x}$, for any $j\ge 4$ we find positive function $\vv_{1}$ with $\vv_{1}(r)\downarrow 0$ as $r\downarrow 0$ such that for $\|\phi\|_{L^k(\mu)}\le 1$,
\beq\label{N04} \E\bigg[\bigg(\int_0^T \big\{\big|\eta_s^r|^2+\|\xi_s^r\|^2\big\}\d s\bigg)^j\bigg]\le \vv_{1}(r)(1+ |\phi(x)|^{2j}),\ \ r\in [0,1], x\in\R^d.\end{equation}
Combining this with \eqref{NPP}, $(H)$ and BDG's inequality, we find a constant $c_1>0$ such that
$$\gg_t^r:= \sup_{s\in [0,t]}|z_s^r|,\ \ t\in [0,T]$$ satisfies
\beq\label{N05}  \E[\gg_t^{j} ] \le \vv_2(r)+ c_1 \int_0^t\Big\{\gg_s^j + \gg_s^{j-1}|\eta_s^r| + \gg_s^{j-2} |\xi_s^r|^2 \Big\}\d s,\ \ t\in [0,T],\end{equation}  where by \eqref{AC''}
$\vv_2(r):= \E[|z_0^r|^j]\to 0$ as $r\to 0$.
Since $st\le s^{\ff n {n-1}}+ t^n$ holds for $s,t\ge 0$ and $n\ge 1$, by taking $n=\ff j 2 $ and $j$ for $j\ge 4$ respectively, we obtain
\beg{align*} &  \int_0^t\Big\{  \gg_s^{j-1}|\eta_s^r| + \gg_s^{j-2} |\xi_s^r|^2 \Big\}\d s\\
&\le   \bigg(\int_0^t |z_s^r|^{2(j-1)} \d s\bigg)^{\ff 1 2}  \bigg(\int_0^t |\eta_s^r|^2 \d s\bigg)^{\ff 1 2} +   \bigg(\int_0^t |z_s^r|^{2(j-2)} \d s\bigg)^{\ff 1 2}  \bigg(\int_0^t |\xi_s^r|^2 \d s\bigg)^{\ff 1 2} \\
&\le   \bigg(\int_0^t |z_s^r|^{2(j-1)} \d s\bigg)^{\ff j{2(j-1)}} +   \bigg(\int_0^t |z_s^r|^{2(j-2)} \d s\bigg)^{\ff j{2(j-2)}}+\aa^r,\end{align*}
where
$$\aa^r:=  \bigg(\int_0^T |\eta_s^r|^2 \d s\bigg)^{\ff j 2} + \bigg(\int_0^T |\xi_s^r|^2 \d s\bigg)^{j}.$$
So, there exists a constant $c_2>0$ such that
\beg{align*}  &c_1 \int_0^t\Big\{  \gg_s^{j-1}|\eta_s^r| + \gg_s^{j-2} |\xi_s^r|^2 \Big\}\d s\\
&\le  c_1 |\gg_t^r|^{\ff{j(j-2)}{2(j-1)}} \bigg(\int_0^t |z_s^r|^j \d s\bigg)^{\ff j{2(j-1)}} + c_1|\gg_t^r|^{\ff{j(j-4)}{2(j-2)}}  \bigg(\int_0^t |z_s^r|^j \d s\bigg)^{\ff j{2(j-2)}}
+ c_1 \aa^r\\
&\le \ff 1 2 |\gg_t^r|^j+ c_2 \int_0^t |\gg_s^r|^j\d s + c_1 \aa^r. \end{align*}
Since \eqref{DFF'} implies $\E[|\gg_t^r|^j]<\infty,$ combining this with \eqref{N04}, \eqref{N05} and applying Gronwall's inequality, we prove \eqref{LA1'} for some positive function $\vv$ with $\vv(r)\downarrow 0$ as $r\downarrow 0$.
 \end{proof}

We are now ready to prove $\aa_i(r)\to 0$ as $r\to 0$ for $i=1,2,3,4$ respectively and hence finish the proof of Theorem \ref{TA2''}.

\ \newline {\bf (a)}  $ \aa_1(r)\to 0$.  As in \eqref{N0N}, by $(H)$ and \eqref{00*} we find  a sequence of positive numbers $s_n\downarrow 0$ as $n\uparrow \infty$ such that
\beq\label{11*}\sup_{s\in [0,T]}  \|\zeta_s(x)-\zeta_s(y)\|^2\le   n |x-y|^{2(k-1)}+ s_n,\ \ n\ge 1,\end{equation}
\beq\label{22*} \beg{split} &\sup_{s\in [0,T]} \big\|D^L b_s^{(1)}(x,\mu)(y)-D^L b_s^{(1)}(x',\nu)(y')\big\|\\
&\le   n \{|x-x'|+|y-y'|+\W_k(\mu,\nu)\}^{\ff 1 {k^*}}+ s_n,\ \ n\ge 1.\end{split}\end{equation}
By \eqref{11*}, Lemma \ref{LM1}, Lemma \ref{LM2} and \eqref{DFF'-}, we find a constant $c_1>0$ such that for any $\phi\in T_{\mu,k}$ with $\|\phi\|_{L^k(\mu)}\le 1$,
\beg{align*} &\E[|F_\phi(r,x)-F_\phi(0,x)|]\\
&\le \E\bigg(\int_0^t |\zeta_s(X_s^{\mu^r,x+r\phi(x)}) - \zeta_s(X_s^{\mu,x}) |^2\cdot |\nn_{\phi(x)}X_s^{\mu^r, x+r\phi(x)}|^2\d s\bigg)^{\ff 1 2}\\
&\quad +\|\zeta\|_\infty \E\bigg(\int_0^t   |\nn_{\phi(x)}X_s^{\mu^r, x+r\phi(x)}-\nn_{\phi(x)}X_s^{\mu, x} |^2\d s\bigg)^{\ff 1 2}\\
&\le \bigg(\E\Big[\sup_{s\in [0,T]} |\nn_{\phi(x)} X_s^{\mu^r,x+r\phi(x)}|^2\Big] \int_0^t \E\big[n|X_s^{\mu^r,x+r\phi(x)}-X_s^{\mu,x}|^{2(k-1)}+ s_n\big]\d s\bigg)^{\ff 1 2} \\
&\quad  +c_1|\phi(x)|  \min\big\{1, \vv(r)(1+|\phi(x)|)\big\} \\
&\le c_1 |\phi(x)| \Big( \ss n   (r+r|\phi(x)|)^{k-1}+\ss{s_n} + \min\big\{1, \vv(r)(1+|\phi(x)|)\big\}\Big)\ \ n\ge 1.\end{align*}
Integrating with respect to $\mu(\d x)$ and  letting first $r\to 0$ then $n\to\infty$, we prove $\aa_1(r)\to 0$ as $r\to 0.$

\ \newline
{\bf (b) }    $\aa_2(r)+\aa_4(r)\to 0$. Let
$$R_\theta := \int_0^\theta \<\zeta_s(X_s^{\mu,x})\nn_{\phi(x)}X_s^{\mu,x}, \d W_s\>,\ \ \theta\in [0,t].$$
By \eqref{DFF'}, we find a constant $c_1>0$ such that
\beq\label{G1} \E[|R_t-R_\theta|]\le c_1 \ss{t-\theta}|\phi(x)|,\ \ \theta\in [0,t], x\in\R^d.\end{equation}
On the other hand, as in \eqref{NL0} and \eqref{NL*}, we find a constant $c_2>0$ such that for $\|\phi\|_{L^k(\mu)}\le 1$,
\beg{align*}&\big|\E\big[f(X_t^{\mu^r, x+r\phi(x)}) - f(X_t^{\mu,x})\big|\F_0\big] \big|= \big|(P_{\theta,t}^{\mu^r} f)(X_\theta^{\mu^r, x+r\phi(x)})- (P_{\theta,t}^\mu f)(X_\theta^{\mu,x})\big|\\
&\le \big|(P_{\theta,t}^{\mu^r} f)(X_\theta^{\mu^r, x+r\phi(x)})- (P_{\theta,t}^{\mu^r} f)(X_\theta^{\mu,x})\big| + \big|(P_{\theta,t}^{\mu^r} f)(X_\theta^{\mu, x})- (P_{\theta,t}^\mu f)(X_\theta^{\mu,x})\big|\\
&\le c_2 \|f\|_\infty\Big[1\land \ff{|X_\theta^{\mu^r, x+r\phi(x)} -X_\theta^{\mu,x}|}{\ss{t-\theta}} + r\Big].\end{align*}
Combining this with \eqref{G1} and Lemma \ref{LM1}, we find   constants $c_3, c_4>0$ such that
\beg{align*} & \big|\E\big[\{f(X_t^{\mu^r, x+r\phi(x)}) - f(X_t^{\mu,x})\} F_\phi(0,x)\big]\big|\\
&\le 2 \|f\|_\infty \E[|R_t-R_\theta|] + \big|\E\big[\E\big(f(X_t^{\mu^r, x+r\phi(x)}) - f(X_t^{\mu,x})\big|\F_\theta\big)R_\theta\big]\big|\\
&\le c_3\|f\|_\infty \bigg\{ \ss{t-\theta}|\phi(x)|+ \Big( r + \ff{\min\{1, r(1+|\phi(x)|)\}}{\ss{t-\theta}}\Big)\big(\E[|R_\theta|^2]\big)^{\ff 1 2} \bigg\}\\
&\le c_4\|f\|_\infty\bigg\{\ss{t-\theta} |\phi(x)| + r|\phi(x)| + \ff{\{n r^{k-1} (1+|\phi(x)|)^{k-1} + s_n\}|\phi(x)|}{\ss{t-\theta}}\bigg\},\end{align*}
where
$$s_n:= \sup_{s>0}\{s\land 1 - n s^{k-1}\}\downarrow 0\ \text{as}\ n\uparrow\infty.$$
Therefore, there exists a constant $c_5>0$ such that
$$\aa_2(r)\le \|f\|_\infty\bigg\{c_5\ss{t-\theta} +\ff{n r^{k-1}+s_n}{\ss{t-\theta}} + r\bigg\},\ \ \theta\in (0,t).$$
By letting first $r\to 0$ then $n\to\infty$ and finally $\theta\to t$, we prove $\aa_2(r)\to 0$ as $r\to 0$.

The proof of $\aa_4(r)\to 0$ is completely similar.

\ \newline
{\bf (c)}    $\aa_3(r)\to 0$.  Write
$$\E\big[|G_\phi(r)-G_\phi(0)|\big] \le \vv_r(\phi)   +\|\zeta\|_\infty \E\big[J_r(X_s^{\mu^r}, X_s^\mu)\big],$$
where
\beg{align*} &  \vv_r(\phi):= \E\bigg[\bigg(\int_0^t |\zeta_s(X_s^r)-\zeta_s(X_s)|^2 \big(\E|\nn_{\phi(X_0)}X_s^{\mu^r}|^{k}\big)^{\ff 2 {k}}\d s \bigg)^{\ff 1 2}\bigg],\\
& J_r(y,z):= \bigg(\int_0^t \big(\E\big[\<D^Lb_s^{(1)}(y_s, \mu_s^r)(X_s^r), \nn_{\phi(X_0)}X_s^{\mu^r}\> \\
&\qquad \qquad \qquad - \<D^Lb_s^{(1)}(z_s, \mu_s)(X_s), \nn_{\phi(X_0)}X_s^{\mu}\>\big]\big)^2 \d s\bigg)^{\ff 1 2},\ \ y,z\in C([0,t];\R^d).\end{align*}
By \eqref{DFF-} for $j=k$,  we obtain
\beq\label{LSTY}  \sup_{r\in [0,1]} \E[|\nn_{\phi(X_0)}X_s^{\mu^r}|^{k}] \le c,\ \ \|\phi\|_{L^k(\mu)}\le 1\end{equation}  for some constant $c>0$, so that by
\eqref{EE0*} and \eqref{11*}, we find   constants $c_1, c_2>0$ such that
$$  \sup_{\|\phi\|_{L^k(\mu)}\le 1}  \vv_r(\phi) \le c_1 \E\bigg[\sup_{s\in [0,t]} n |X_s^r-X_s|^{k-1}+ s_n\bigg] \le c_2 n r^{k-1} + c_1s_n,\ \ n\ge 1.$$
Then  $ \sup_{\|\phi\|_{L^k(\mu)}\le 1} \vv_r(\phi) \to 0$ as $r\to 0$. It remains to prove
\beq\label{NMM} \lim_{r\downarrow 0} \sup_{\|\phi\|_{L^k(\mu)}\le 1}  \E[J_r(X^r,X)]= 0.\end{equation}
By $(H)$, \eqref{EE0*}, \eqref{LN2}, Lemma \ref{LM2}, \eqref{22*} and \eqref{LSTY}, we find constants $c_3,c_4, c_5>0$ and positive function $\tt\vv(\cdot)$ on $[0,1]$ with $\tt\vv(r)\to 0$ as $r\to 0$, such that when $\|\phi\|_{L^k(\mu)}\le 1$,
\beg{align*} &\E\big[\big|\<D^Lb_s^{(1)}(y_s, \mu_s^r)(X_s^r), \nn_{\phi(X_0)}X_s^{\mu^r}\> - \<D^Lb_s^{(1)}(z_s, \mu_s)(X_s), \nn_{\phi(X_0)}X_s^{\mu}\>\big|\big]\\
&\le c_3\big(\E\big[ \big|\nn_{\phi(X_0)} X_s^{\mu^r}-  \nn_{\phi(X_0)}X_s^{\mu}\big|^k\big]\big)^{\ff 1 k}  \\
&\quad + \Big(\E[|\nn_{\phi(X_0)}X_s|^k\big]\big)^{\ff 1 k}
\big(\E[|D^Lb_s^{(1)}(z_s, \mu_s)(X_s)-D^Lb_s^{(1)}(y_s, \mu_s^r)(X_s^r)|^{k^*} ]\big)^{\ff 1 {k^*}} \\
& \le  \tt\vv(r)  + c_4 \big(\E[ n^{k^*} \{|z_s-y_s|+  |X_s^r-X_s|+ r\} +s_n^{k^*}]\big)^{\ff 1 {k^*}}  \\
&\le \tt\vv(r)  + c_5\big\{n |z_s-y_s|^{\ff 1 {k^*}}+ n r^{\ff 1 {k^*}} + s_n\big\},\ \ n\ge 1. \end{align*}
Combining this with \eqref{EE0*} we find a constant $c_6>0$ such that
$$ \sup_{\|\phi\|_{L^k(\mu)}\le 1}  \E[J_r(X^r,X)]\le c_6\big\{\tt\vv(r)+ n r^{\ff 1{k^*}} + s_n\big\},\ \ n\ge 1.$$
By letting first $r\to 0$ then $n\to\infty$ we prove \eqref{NMM}.

  \paragraph{Acknowledgement.} The author  would like to thank the referee for helpful comments and corrections.

\end{document}